\documentclass{gtart}

\def\ifplaintex{\expandafter\ifx\csname documentclass\endcsname\relax}

\def\gtp{{\mathsurround=0pt\it $\cal G\mskip-2mu$eometry \&\ 
$\cal T\!\!$opology $\cal P\!$ublications}}  

\def\recd{{\small Received:\qua\receiveddate\ifx\reviseddate\relax
\else\qquad Revised:\qua\reviseddate\fi\par}} 

\def\lognumber#1{\def\thelognumber{#1}}
\def\volumenumber#1{\def\thevolumenumber{#1}}
\def\volumeyear#1{\def\thevolumeyear{#1}}
\def\papernumber#1{\def\thepapernumber{#1}}
\def\pagenumbers#1#2{\def\startpage{#1}\def\finishpage{#2}}
\def\published#1{\def\publishdate{#1}}

\def\received#1{\def\receiveddate{#1}}
\def\revised#1{\def\reviseddate{#1}}
\def\accepted#1{\def\accepteddate{#1}}

\def\asciiaddress#1{\def\theasciiaddress{#1}}
\def\asciiemail#1{\def\theasciiemail{#1}}

\long\def\asciiabstract#1{\long\def\theasciiabstract{#1}}

\let\\\par\let\thelognumber\relax\let\thevolumenumber\relax
\let\thepapernumber\relax\let\thevolumeyear\relax\let\startpage\relax
\let\finishpage\relax\let\publishdate\relax\let\receiveddate\relax
\let\reviseddate\relax\let\accepteddate\relax\let\theasciititle\relax
\let\theasciiauthors\relax\let\theasciiaddress\relax
\let\theasciiabstract\relax

\let\theasciiemail\relax

\ifplaintex
\font\logobig=cmssbx10 scaled 3836
\font\logomed=cmssbx10 scaled 2557
\else
\font\logobig=cmssbx10 scaled 4200
\font\logomed=cmssbx10 scaled 2800
\fi

\long\def\makeagttitle{   
\count0=\startpage
\agt\hfill      
\hbox to 45truept{\vbox to 0pt{\vglue -13truept{\logomed A\kern -.37em{\logobig 
T}\kern -.38em G}\vss}\hss}
\break
{\small Volume \thevolumenumber\ (\thevolumeyear)
\startpage--\finishpage\nl
Published: \publishdate}

\vglue .25truein

{\parskip=0pt\leftskip 0pt plus
1fil\def\\{\par\smallskip}{\Large\bf\thetitle}\par\medskip} \vglue
0.05truein

{\parskip=0pt\leftskip 0pt plus 1fil\def\\{\par}{\sc\theauthors}
\par\medskip}%
 
\vglue 0.03truein 

{\small\leftskip 25truept\rightskip 25truept{\bf Abstract}\stdspace\theabstract

{\bf AMS Classification}\stdspace\theprimaryclass
\ifx\thesecondaryclass\relax\else; \thesecondaryclass\fi\par
{\bf Keywords}\stdspace \thekeywords\par}\vglue 7truept

}   

\ifplaintex
\font\phead=cmsl9 scaled 950
\font\pnum=cmbx10 scaled 913
\font\pfoot=cmsl9 scaled 950
\headline{\vbox to 0pt{\vskip -4.5mm\line{\small\phead\ifnum
\count0=\startpage ISSN 1472-2739 (on-line) 1472-2747 (printed)
\hfill {\pnum\folio}\else\ifodd\count0\def\\{ }%
\ifx\theshorttitle\relax\thetitle\else\theshorttitle\fi\hfill{\pnum\folio}
\else\def\\{ and }{\pnum\folio}\hfill\ifx\theshortauthors\relax\theauthors
\else\theshortauthors\fi\fi\fi}\vss}}
\footline{\vbox to 0pt{\vglue 0mm\line{\small\pfoot\ifnum\count0=\startpage
\copyright\ \gtp\hfill\else
\agt, Volume \thevolumenumber\ (\thevolumeyear)\hfill\fi}\vss}}
\else
\font=cmsl9 scaled 1050
\font\lnum=cmbx10 
\font=cmsl9 scaled 1050
\makeatletter
\def\@oddhead{{\small\lhead\ifnum\count0=\startpage ISSN 1472-2739 
(on-line) 1472-2747 (printed)\hfill {\lnum\number\count0}\else\ifodd\count0
\def\\{ }\ifx\theshorttitle\relax \thetitle \else\theshorttitle\fi\hfill
{\lnum\number\count0}\else\def\\{ and }{\lnum\number\count0}
\hfill\ifx\theshortauthors\relax 
\theauthors\else\theshortauthors\fi\fi\fi}}\def\@evenhead{\@oddhead}
\def\@oddfoot{\small\lfoot\ifnum\count0=\startpage\copyright\ \gtp\hfill\else
\agt, Volume \thevolumenumber\ (\thevolumeyear)\hfill\fi}
\def\@evenfoot{\@oddfoot}
\makeatother
\fi
\let\maketitlepage\makeagttitle
\let\makeshorttitle\maketitlepage
\let\maketitle\maketitlepage

\newwrite\gtoutfile
\long\gdef\makeheadfile{  
{\def\\{, }\def\s{ }
\immediate\openout\gtoutfile head.xxx
\immediate\write\gtoutfile{To: math@arxiv.org}
\immediate\write\gtoutfile{Subject: put OR rep NNNNN:ppppp}
\immediate\write\gtoutfile{--text follows this line--}
\immediate\write\gtoutfile{Proxy-for: \ifx\theasciiauthors\relax
\theauthors\else\theasciiauthors\fi\s<\ifx\theasciiemail\relax\theemail\else\theasciiemail\fi>}
\immediate\write\gtoutfile{\noexpand\\}
\immediate\write\gtoutfile{Authors: \ifx\theasciiauthors\relax
\theauthors\else\theasciiauthors\fi}
{\def\\{ }\immediate\write\gtoutfile{Title: \ifx\theasciititle\relax
\thetitle\else\theasciititle\fi}}
\immediate\write\gtoutfile{Subj-class: GT or SG, GR etc}
\immediate\write\gtoutfile{MSC-class: \theprimaryclass\ifx\thesecondaryclass\relax\else, \thesecondaryclass\fi}
\immediate\write\gtoutfile{Journal-ref: Algebr. Geom. Topol. \thevolumenumber\s
(\thevolumeyear) \startpage-\finishpage}
\immediate\write\gtoutfile{Comments: Published by Algebraic and
Geometric Topology at}
\immediate\write\gtoutfile{\s\s\s  http://www.maths.warwick.ac.uk/agt/AGTVol\thevolumenumber/agt-\thevolumenumber-\thepapernumber.abs.html}
\immediate\write\gtoutfile{\noexpand\\}
\immediate\write\gtoutfile{}
\ifx\theasciiabstract\relax
\immediate\write\gtoutfile{\theabstract}\else
\immediate\write\gtoutfile{\theasciiabstract}\fi
\immediate\write\gtoutfile{}
\immediate\write\gtoutfile{\noexpand\\}
\immediate\write\gtoutfile{}
\immediate\closeout\gtoutfile}}  

\def\maketitlepage{\makeagttitle\makeheadfile}
\let\makeshorttitle\maketitlepage
\let\maketitle\maketitlepage

\lognumber{11}
\volumenumber{3}
\volumeyear{2003}
\papernumber{11}
\published{13 March 2003}
\pagenumbers{287}{334}
\received{18 August 2002}
\revised{11 February 2003}
\accepted{11 March 2003}

\usepackage{amsmath,amssymb}
\usepackage[matrix,arrow,curve]{xy}

\numberwithin{equation}{section}
\newtheorem{theorem}[equation]{Theorem}
\newtheorem{lemma}[equation]{Lemma}
\newtheorem{proposition}[equation]{Proposition}
\newtheorem{corollary}[equation]{Corollary}

\theoremstyle{definition}
\newtheorem{definition}[equation]{Definition}
\newtheorem{remark}[equation]{Remark}
\newtheorem{example}[equation]{Example}

\DeclareMathOperator{\Ho}{Ho}
\DeclareMathOperator{\Ev}{Ev}
\DeclareMathOperator{\map}{map}

\newcommand{\mI}{{\mathbb I}}
\newcommand{\mN}{{\mathbb N}}
\newcommand{\mS}{{\mathbb S}}
\newcommand{\mT}{{\mathbb T}}

\newcommand{\mP}{{\mathbb P}}
\newcommand{\mU}{{\mathbb U}}
\newcommand{\mZ}{{\mathbb Z}}
\newcommand{\mIC}{{\mathbb I}_{\mathcal C}}
\newcommand{\mID}{{\mathbb I}_{\mathcal D}}

\newcommand{\cA}{{\mathcal A}}

\newcommand{\cC}{{\mathcal C}}
\newcommand{\cD}{{\mathcal D}}

\newcommand{\cF}{{\mathcal F}}
\newcommand{\cI}{{\mathcal I}}
\newcommand{\cM}{{\mathcal M}}
\newcommand{\cO}{{\mathcal O}}
\newcommand{\cR}{{\mathcal R}}
\newcommand{\cT}{{\mathcal T}}
\newcommand{\cS}{{\mathcal S}}

\newcommand{\cW}{{\mathcal W}}

\newcommand{\Hom}{\mathrm{Hom}}
\newcommand{\alg}{{\mbox{-}\mathcal A}lg}
\newcommand{\mon}{{\mbox{-}\mathcal M}onoid}
\newcommand{\boxprod}{\mathbin{\square }}
\newcommand{\cat}{{\mbox{-}\mathcal C}at}
\newcommand{\ch}{ch^+}
\newcommand{\dga}{{\mathcal DGA}}
\newcommand{\dgr}{{\mathcal DGR}}
\newcommand{\graph}{\mbox{-}{\mathcal G}raph}
\newcommand{\iso}{\cong}
\newcommand{\sab}{s{\mathcal A}b}
\newcommand{\sm}{\wedge}
\newcommand{\spec}{Sp^{\Sigma}}
\newcommand{\sring}{s{\mathcal R}}
\newcommand{\tensor}{\otimes}
\newcommand{\tvarphi}{\widetilde{\varphi}}
\newcommand{\Gam}{\Gamma}
\newcommand{\Id}{\text{Id}}

\newcommand{\Modr}{{\mathcal M}od\mbox{-}}
\renewcommand{\to}{\longrightarrow}

\begin{document}

\title{Equivalences of monoidal model categories}

\authors{Stefan Schwede\\Brooke Shipley}
\addresses{SFB 478 Geometrische Strukturen in der Mathematik\\ 
Westf\"alische Wilhelms-Universit\"at M\"unster, Germany \medskip\\
Department of Mathematics, Purdue University\\W. Lafayette, IN 47907, USA}
\asciiaddress{SFB 478 Geometrische Strukturen in der Mathematik\\ 
Westfaelische Wilhelms-Universitaet Muenster, Germany\\and\\
Department of Mathematics, Purdue University\\W. Lafayette, IN 47907, USA}
\email{sschwede@math.uni-muenster.de}
\secondemail{bshipley@math.purdue.edu}
\asciiemail{sschwede@math.uni-muenster.de, bshipley@math.purdue.edu}

\begin{abstract}
We construct Quillen equivalences between the model
categories of monoids (rings), modules and algebras over two Quillen
equivalent model categories under certain conditions.  This is
a continuation of our earlier work where we established model
categories of monoids, modules and algebras~\cite{ss}.  As an application
we extend the Dold-Kan equivalence to show that the model
categories of simplicial rings, modules and algebras are Quillen
equivalent to the associated model categories of connected differential
graded rings, modules and algebras.  
We also show that our classification results from~\cite{stable}
concerning stable model categories translate to any one 
of the known symmetric monoidal model categories of spectra.
\end{abstract}

\asciiabstract{We construct Quillen equivalences between the model
categories of monoids (rings), modules and algebras over two Quillen
equivalent model categories under certain conditions.  This is a
continuation of our earlier work where we established model categories
of monoids, modules and algebras [Algebras and modules in monoidal
model categories, Proc. London Math. Soc. 80 (2000), 491-511].  As an
application we extend the Dold-Kan equivalence to show that the model
categories of simplicial rings, modules and algebras are Quillen
equivalent to the associated model categories of connected
differential graded rings, modules and algebras.  We also show that
our classification results from [Stable model categories are
categories of modules, Topology, 42 (2003) 103-153] concerning stable
model categories translate to any one of the known symmetric monoidal
model categories of spectra.}

\primaryclass{55U35}
\secondaryclass{18D10, 55P43, 55P62}
\keywords{Model category, monoidal category, Dold-Kan equivalence, spectra}

\maketitle

\section{Introduction}
This paper is a sequel to~\cite{ss} where we studied sufficient conditions
for extending Quillen model category structures to the associated categories
of monoids (rings), modules and algebras over a monoidal model category. 
Here we consider functors between such categories.  
We give sufficient conditions for extending Quillen
equivalences of two monoidal model categories to Quillen equivalences
on the associated categories of monoids, modules and algebras.  
This is relatively easy when the initial Quillen equivalence is via an
adjoint pair of functors which induce adjoint functors on the 
categories of monoids; see for example~\cite[\S 13, 16]{mmss} 
and~\cite[5.1]{sch}. We refer to this situation as a
{\em strong monoidal Quillen equivalence}, 
see Definition \ref{def-WMQ equivalence}.

However, in the important motivating example of chain
complexes and simplicial abelian groups, only something weaker holds:
the right adjoint has a monoidal structure,
but the left adjoint only has a lax {\em co}monoidal 
(also referred to as op-lax monoidal) structure which
is a weak equivalence on cofibrant objects.
We refer to this situation as a {\em weak monoidal Quillen equivalence}.
Our general result about monoidal equivalences, 
Theorem~\ref{thm-WMQE gives equivalent algebras}, 
works under this weaker assumption.
Also, in Proposition~\ref{prop-generating criterion} 
we give a sufficient criterion for showing that an adjoint pair 
is a weak monoidal Quillen equivalence.

Our motivating example is the Dold-Kan equivalence of chain
complexes and simplicial abelian groups.  The normalization
functor $N \co \sab \to \ch$ is monoidal with respect to the
graded tensor product of chains, the levelwise tensor product of
simplicial abelian groups and the transformation 
known as the {\em shuffle map}; the inverse equivalence $\Gam\co \ch \to \sab$ 
also has a monoidal structure (coming from the {\em Alexander-Whitney map}).

The natural isomorphism $N\Gam\iso\Id_{\ch}$ is monoidal with respect
to the shuffle and Alexander-Whitney maps. 
This implies that the algebra valued functor $\Gam$ 
embeds connective differential graded rings
as a full, reflexive subcategory of the category of simplicial rings,
see Proposition \ref{prop-full embedding}. 
However, the other natural isomorphism $\Gam N\iso\Id_{\sab}$
cannot be chosen in a monoidal fashion.
Hence these functors do not induce inverse, or even adjoint functors 
on the categories of algebras.  
One of the main points of this paper is to show that nevertheless,
the homotopy categories of simplicial rings and 
connective differential graded rings are equivalent, via
a Quillen equivalence on the level of model categories.
This Quillen equivalence should be well known 
but does not seem to be in the literature.
A similar equivalence, between reduced rational simplicial Lie algebras 
and reduced rational differential graded Lie algebras, 
was part of Quillen's work on rational homotopy theory~\cite[I.4]{Q2} 
which originally motivated the definition of model categories.

In the following theorem we use the word {\em connective} for 
{\em non-negatively graded} or more precisely {\em $\mN$-graded}
objects such as chain complexes or algebras.

\begin{theorem}\label{thm-sab-ch}       

\begin{enumerate}
\item\label{i} For a connective differential graded ring $R$, 
there is a Quillen equivalence between the categories 
of connective $R$-modules and simplicial modules over 
the simplicial ring $\Gam R$, 
\[\Modr R \ \simeq_Q \ \Modr \Gam R \ . \] 

\item\label{ii} For a simplicial ring $A$ there is a Quillen 
equivalence between the categories of 
connective differential graded $NA$-modules and simplicial modules over $A$,
\[ \Modr NA \ \simeq_Q \ \Modr A \ . \]

\item\label{iii} For a commutative ring $k$,
there is a Quillen equivalence between the categories of connective 
differential graded $k$-algebras and simplicial $k$-algebras,
\[ \dga_k \ \simeq_Q \ k\alg \ . \] 

\item\label{iv} For $A$ a simplicial commutative ring, there is
a Quillen equivalence between  the categories of
connective differential graded $NA$-algebras and simplicial $A$-algebras,
\[NA\alg \ \simeq_Q \ A\alg \ . \] 
\end{enumerate}
\end{theorem}       

The special case for $k=\mZ$ in part \eqref{iii} of the previous theorem
in particular says that the model categories of connective 
differential graded rings and simplicial rings are Quillen equivalent.
Part \eqref{iii} is a special case of part \eqref{iv}
for $A$ a constant commutative simplicial ring. 
The proof of Theorem~\ref{thm-sab-ch} is an application
of the more general Theorem \ref{thm-WMQE gives equivalent algebras}.
Parts \eqref{ii} and \eqref{iii} of Theorem~\ref{thm-sab-ch}
are established in Section \ref{N as right adjoint},
part \eqref{i} is shown in Section \ref{Gamma as right adjoint},
and part \eqref{iv} is completed in~\ref{algebras over com sRing}.

In part \eqref{i} of Theorem \ref{thm-sab-ch}, 
the right adjoint of the Quillen equivalence is induced 
by the functor $\Gamma$ from connective chain complexes 
to simplicial abelian groups 
which is inverse to  normalized chain complex functor $N$.
However, the left adjoint is in general {\em not} given by
the normalized chain complex on underlying simplicial abelian groups.
In parts \eqref{ii} to \eqref{iv}, the right adjoint of the
Quillen equivalence is always induced by the normalized chain complex
functor. However, the left adjoint is in general {\em not} given by
the functor $\Gamma$ on underlying chain complexes.
We discuss the various left adjoints in Section \ref{sub-left adjoint}.

\medskip

Notice that we do not compare the categories 
of {\em commutative} simplicial rings and {\em commutative}
differential graded rings. The normalization functor is symmetric monoidal
with respect to the shuffle map. Hence it takes commutative simplicial
rings to commutative (in the graded sense) differential graded rings.
But the Alexander-Whitney map is not symmetric, and so $\Gamma$ does not
induce a functor backwards.
In characteristic zero, i.e., for algebras $k$ over the rational numbers,
there is a model structure on commutative differential graded rings 
with underlying fibrations and trivial fibrations~\cite{BG, stanley}. 
Moreover, the normalized chain complex functor is then the right adjoint
of a Quillen equivalence between commutative simplicial $k$-algebras and
connective differential graded $k$-algebras; indeed, as Quillen
indicates on p.~223 of~\cite{Q2}, a similar method as for rational
Lie algebras works for rational commutative algebras.

Without a characteristic zero assumption, 
not every commutative differential graded ring 
is quasi-isomorphic to the normalization of a commutative simplicial ring:
if $A$ is a commutative simplicial ring, then every element $x$ 
of odd degree in the homology algebra $H_*(NA)$ satisfies $x^2=0$; 
but in a general commutative differential graded algebra 
we can only expect the relation $2\cdot x^2=0$. More generally,
the homology algebra  $H_*(NA)$, for $A$ a commutative simplicial ring,
has divided power~\cite{cartan} 
and other operations~\cite{Dwyer-operations} which need not be supported 
by a general commutative differential graded algebra.
Moreover, in general the forgetful functor from differential graded algebras 
to chain complexes does not create a model structure and there is no 
homotopically meaningful way to go from differential graded 
to simplicial algebras in a way that preserves commutativity.
While the normalization functor on commutative algebras 
still has a left adjoint, it is not clear if that adjoint preserves 
enough weak equivalences and whether it admits a derived functor.

In arbitrary characteristic, one should consider the categories 
of $E_{\infty}$-algebras instead of the commutative algebras.  
Mandell~\cite[1.2]{mandell} establishes a Quillen equivalence,
in any characteristic, between $E_{\infty}$-simplicial algebras and connective
$E_{\infty}$-differential graded algebras.
The symmetry properties of the Dold-Kan equivalence were also studied
by Richter \cite{richter-DoldKan}; she has shown that for every 
differential graded algebra $R$ which is commutative (in the graded sense),
the simplicial ring $\Gam A$ admits a natural $E_{\infty}$-multiplication.

\medskip

One of the reasons we became interested in generalizing the Dold-Kan 
equivalence is because it is the basis for one out of the four steps 
of a zig-zag of weak monoidal Quillen equivalences between 
$H\mZ$-modules, and $\mZ$-graded chain complexes; see~\cite{sh-Q}. 
Theorem~\ref{thm-WMQE gives equivalent algebras}
then applies to each of these four steps to produce
Quillen equivalences between $H\mZ$-algebras and $\mZ$-graded
differential graded algebras and between the associated module 
categories~\cite[1.1]{sh-Q}.  
These equivalences then provide an algebraic model for 
any rational stable model category with a set of small generators.
These rational algebraic models,~\cite[1.2]{sh-Q},
are really unfinished business from~\cite{stable} and even appeared
in various preprint versions.  These models 
are used as stepping stones in~\cite{S-QT} and~\cite{GS} to
form explicit, small algebraic models for the categories of rational
$T^n$-equivariant spectra for $T^n$ the $n$-dimensional torus.

Another motivation for this general approach to monoidal Quillen equivalences
is the extension of our work in~\cite{stable} where we characterize
stable model categories with a set of generators as those model categories
which are Quillen equivalent to modules over a symmetric ring spectrum
with many objects ($\spec$-category).   Here we show this characterization 
can be translated to any of the other symmetric monoidal categories of spectra.
Quillen equivalences of monoids, modules and algebras 
for these categories of spectra were
considered in~\cite{mmss} and~\cite{sch}, but ring spectra with
many objects (or enriched categories) were not considered. 
Using the Quillen equivalences between modules over ring spectra with
many objects (or enriched categories) 
over the other known highly structured categories of spectra 
established in Section \ref{sec-spectra}, the 
characterization of~\cite[Thm 3.3.3]{stable}
can be translated to any other setting: 

\begin{corollary}\label{cor-spectra}
Let $\cC$ be any of the monoidal model categories 
of symmetric spectra~\cite{hss, mmss}
(over simplicial sets or topological spaces), orthogonal spectra~\cite{mmss}, 
$\cW$-spaces~\cite{mmss}, simplicial functors~\cite{lydakis} or 
$S$-modules~\cite{ekmm}.
Then any cofibrantly generated, proper, simplicial, stable model category
with a set of small generators is Quillen equivalent to modules over a 
$\cC$-category with one object for each generator.
\end{corollary}

{\bf Organization}\qua  
In Section~\ref{sec-sAb and ch} we motivate our general results by
considering the special case of chain complexes and simplicial
abelian groups.  We then turn to the general case and
state sufficient conditions for extending
Quillen equivalences to monoids, modules and algebras in 
Theorem~\ref{thm-WMQE gives equivalent algebras}.  
Section \ref{subsec-criterion} gives a criterion for a Quillen functor pair 
to be weakly monoidal: a sufficient condition
is that one of the unit objects detects weak equivalences.
In Section~\ref{sec-sAb and ch revisited} we return to
chain complexes and simplicial abelian groups and deduce
Theorem~\ref{thm-sab-ch} from the general result. 
Section~\ref{sec-proofs} contains the proof of the main theorem 
and the proof of the criterion for a Quillen functor pair 
to be weakly monoidal. In Section~\ref{app} we consider rings 
with many objects (enriched categories) and their modules.  In 
Theorem~\ref{thm-WMQE many generators} we extend the Quillen equivalences to
modules over these enriched categories.  In Section~\ref{sec-spectra}
we show that these general statements 
apply to the various symmetric monoidal categories of spectra and
deduce Corollary~\ref{cor-spectra}.  
Throughout this paper, modules over a ring, algebra, category, etc,
are always {\em right} modules.

\medskip
{\bf Acknowledgments}\qua We would like to thank Mike Mandell 
for several helpful conversations.
The second author was partially supported by an NSF Grant.

\section{Chain complexes and simplicial abelian groups}\label{sec-sAb and ch}

As motivation for our general result, in this section we begin the 
comparison of the categories of differential graded rings 
and simplicial rings.  We recall the normalized chain complex functor
$N$, its inverse $\Gamma$ and the shuffle and Alexander-Whitney maps.  
We then consider the monoidal properties of the adjunction unit
$\eta \co \Id \to \Gamma N$ and counit  $\epsilon \co N\Gamma \to \Id$.
As mentioned in the introduction, $\epsilon$ is monoidally better behaved 
than $\eta$.  This motivates developing our general result in 
Section~\ref{sec-main} which does not require monoidal adjunctions.
In Section~\ref{sec-sAb and ch revisited} 
we then revisit this specific example and prove Theorem~\ref{thm-sab-ch}.

\subsection{Normalized chain complexes}
The {\em (ordinary) chain complex} $CA$ of a simplicial abelian group $A$
is defined by $(CA)_n=A_n$ with differential the alternating sum
of the face maps,
\[ d \ = \ \sum_{i=0}^n (-1)^id_i \ : \  (CA)_n \ \to \ (CA)_{n-1} \ . \]
The chain complex $CA$ has a natural subcomplex $DA$, the complex of 
{\em degenerate simplices}; by definition, $(DA)_n$ is the 
subgroup of $A_n$ generated by all degenerate simplices.
The {\em normalized chain complex} $NA$ is the quotient complex of $CA$
by the degenerate simplices,
\begin{equation}\label{def-normalized}
NA \ = \ CA/DA \ . \end{equation}
The degenerate complex $DA$ is acyclic, so the projection $CA\to NA$
is a quasi-isomorphism. 

The complex of degenerate simplices has a natural complement,
sometimes called the {\em Moore complex}. 
The $n$-th chain group of this subcomplex is the intersection 
of the kernels of all face maps, except the 0th one,
and the differential in the subcomplex is thus given by the remaining 
face map $d_0$. The chain complex $CA$ is the internal direct sum
of the degeneracy complex $DA$ and the Moore complex.
In particular, the Moore subcomplex is naturally isomorphic 
to the normalized chain complex $NA$; in this paper, we do not use the Moore
complex.

The normalization functor
\[ N \ : \ \sab \ \to \  \ch \]  
from simplicial abelian groups to non-negatively graded chain complexes 
is an equivalence of categories~\cite[Thm.~1.9]{dold}.
The value of the inverse $\Gam\co \ch\to\sab$ on a complex $C$ can be defined by
\begin{equation}\label{eq-def Gamma}
(\Gam C)_n \ = \ \ch(N\Delta^n,C) \ , \end{equation}
where $N\Delta^n$ is short for the normalized chain complex
of the simplicial abelian group freely generated by the standard $n$-simplex.
The simplicial structure maps in $\Gam C$ are induced from the cosimplicial
structure of $\Delta^n$ as $n$ varies through the simplicial 
category $\Delta$. A natural isomorphism 
\[ \eta_A \ : \ A \ \to \ \Gam NA \]
is defined in simplicial dimension $n$ by
\[ A_n \ni a \ \longmapsto \ (N\bar a\co N\Delta^n\to NA) \in  (\Gam NA)_n  \]
where $\bar a\co \Delta^n\to A$ is the unique morphism of simplicial sets which
sends the generating $n$-simplex of $\Delta^n$ to $a\in A_n$.
The other natural isomorphism $\epsilon_C\co N\Gam C\to C$ is uniquely determined
by the property 
\[ \Gam(\epsilon_C) \ = \ \eta_{\Gam C}^{-1} \ : \ \Gam N\Gam C\to \Gam C\ . \]

\subsection{Tensor products}\label{sub-Tensor products} 
The tensor product of two connective chain complexes $C$ and $D$ 
is defined by
\begin{equation}\label{eq-tensor of complexes} 
(C\tensor D)_n \ = \ \bigoplus_{p+q=n} \ C_p\tensor D_q \end{equation}
with differential given on homogeneous elements
by 
\[ d(x\tensor y) \ = \ dx\tensor y \, + \, (-1)^{|x|}x\tensor dy \ . \]
The tensor product of simplicial abelian groups is defined dimensionwise.

Both tensor products are symmetric monoidal. 
The respective unit object is the free abelian group of rank one,
viewed either as a complex concentrated in dimension zero
or a constant simplicial abelian group. The associativity and unit
isomorphisms are obvious enough that we do not specify them, similarly
the commutativity isomorphism for simplicial abelian groups.
The commutativity isomorphism for complexes involves a sign, i.e.,
\[ \tau_{C,D} \ : \ C\tensor D \ \to \ D\tensor C \]
is given on homogeneous elements
by $\tau(x\tensor y)= (-1)^{|x||y|}y\tensor x$.

The unit objects are preserved under the normalization functor and its inverse.
However, the two tensor products for chain complexes and
simplicial abelian groups are different in an essential way, i.e., the
equivalence of categories given by normalization does {\em not}
take one tensor product to the other. Another way of saying this is
that if we use the normalization functor and its inverse 
to transport the tensor product of simplicial abelian groups 
to the category of connective chain complexes, we obtain
a second monoidal product (sometimes called the {\em shuffle product of
complexes}) which is non-isomorphic, and significantly bigger than,
the tensor product \eqref{eq-tensor of complexes}.
Another difference is that the tensor product \eqref{eq-tensor of complexes}
makes perfect sense for $\mZ$-graded chain complexes, whereas the shuffle
product cannot be extended to a monoidal structure 
on $\mZ$-graded chain complexes in any natural way.
The difference between these two tensor products is responsible for the
fact that the categories of simplicial rings and of connective differential
graded rings are not equivalent.

\subsection{The shuffle and Alexander-Whitney maps}
Even though the tensor products of chain complexes and
simplicial modules do not coincide under normalization, 
they can be related in various ways.
The {\em shuffle map}
\begin{equation}\label{def-shuffle}
\nabla \ : \ CA \tensor CB \ \to \ C(A\tensor B) \ , 
\end{equation}
was introduced by Eilenberg and Mac Lane~\cite[(5.3)]{EM-H(pi)}, 
see also~\cite[VIII 8.8]{ML-homology} or \cite[29.7]{may}. 
For simplicial abelian groups $A$ and $B$ and simplices $a\in A_p$ 
and $b\in B_q$, the image 
$\nabla(a\tensor b)\in C_{p+q}(A\tensor B)=A_{p+q}\tensor B_{p+q}$ 
is given by
\[ \nabla(a\tensor b) \ = \ \sum_{(\mu,\nu)} \, \text{sign}(\mu,\nu) \cdot
s_{\nu}a \tensor s_{\mu}b \ ; \] 
here the sum is taken over all {\em $(p,q)$-shuffles}, i.e.,
permutations of the set $\{0,\dots,p+q-1\}$ which leave the first
$p$ elements and the last $q$ elements in their natural order.
Such a $(p,q)$-shuffle is of the form 
$(\mu,\nu)=(\mu_1,\dots,\mu_p,\nu_1,\dots,\nu_q)$ with
$\mu_1<\mu_2<\cdots<\mu_p$ and $\nu_1<\cdots<\nu_q$, and the
associated degeneracy operators are given by
\[s_{\mu}b= s_{\mu_p}\cdots s_{\mu_1}b \text{\quad respectively \quad } 
s_{\nu}a= s_{\nu_q}\cdots s_{\nu_1}a \ .  \] 

The shuffle map is a {\em lax monoidal transformation},
i.e., it is appropriately unital and associative, 
see~\cite[Thm.\ 5.2]{EM-H(pi)} or \cite[29.9]{may}.
The unit map is the unique chain map  $\eta \co \mZ[0] \to C(\mZ)$ 
which is the identity in dimension 0.

The {\em Alexander-Whitney map}~\cite[VIII 8.5]{ML-homology}, \cite[29.7]{may}
\begin{equation}\label{def-AW} 
AW \ : \ C(A\tensor B) \ \to \ CA \tensor CB
\end{equation}
goes in the direction opposite to the shuffle map;
it is defined for a tensor product of $n$-simplices
$a\in A_n$ and $b\in B_n$ by
\[ AW(a\tensor b) \ = \ \bigoplus_{p+q=n} \,
\tilde d^{p}a \tensor d_0^{q}b \ . \] 
Here the `front face' $\tilde d^{p}\co A_{p+q}\to A_p$
and the `back face' $d_0^{q}\co B_{p+q}\to B_q$
are induced by the injective monotone maps 
$\tilde\delta^{p}\co [p]\to [p+q]$ 
and $\delta_0^{q}\co [q]\to [p+q]$
defined by $\tilde\delta^{p}(i)=i$ 
and $\delta_0^{q}(i)=p+i$.
The Alexander-Whitney map is a {\em lax comonoidal transformation},
(also referred to as an {\em op-lax monoidal transformation}) 
i.e., it is appropriately unital and associative.

Both the shuffle and the Alexander-Whitney map preserve the subcomplexes
of degenerate simplices, 
compare~\cite[Lemma 5.3]{EM-H(pi)} or \cite[29.8, 29.9]{may}.
Hence both maps factor over normalized chain complexes and induce maps
\begin{equation}\label{def-restricted shuffle} 
\nabla \co NA \tensor NB \ \to \ N(A\tensor B) \ \text{\ and\ }\
AW \co N(A\tensor B) \ \to \ NA \tensor NB 
\end{equation}
for which we use the same names.
These restricted maps are again lax monoidal, respectively lax comonoidal, 
and the restricted unit maps are now isomorphisms $\mZ[0] \iso N(\mZ)$.

Moreover, the composite map $AW\circ\nabla\co CA \tensor CB\to CA \tensor CB$
differs from the identity only by degenerate simplices.
Hence on the level of normalized complexes, the composite
\begin{equation}\label{AW nabla=Id}
  NA \tensor NB \ \xrightarrow{\nabla} \ N(A\tensor B) \ 
\xrightarrow{AW} \ NA \tensor NB \end{equation} 
is the identity transformation. The composite of shuffle and Alexander-Whitney
maps in the other order are naturally chain homotopic to the 
identity transformation.
In particular, the shuffle map \eqref{def-shuffle}, the
Alexander Whitney map \eqref{def-AW}  and their
normalized versions \eqref{def-restricted shuffle} 
are all quasi-isomorphisms of chain complexes.

The shuffle map is also {\em symmetric} 
in the sense that for all simplicial abelian groups $A$ and $B$, 
the following square commutes
\[\xymatrix{ CA \tensor CB \ar[d]_{\nabla} \ar[r]^{\tau} &
 CB \tensor CA \ar[d]^{\nabla}  \\
C(A\tensor B) \ar[r]_{C(\tau)}& C(B\tensor A) }\]
where $\tau$ denotes the symmetry isomorphism of the tensor products
of either simplicial abelian groups or chain complexes.
The normalized version \eqref{def-restricted shuffle} of the shuffle
map is symmetric as well. However, the Alexander-Whitney map 
is {\em not} symmetric, nor is its normalized version.

We can turn the comonoidal structure on the normalization functor
given by the Alexander-Whitney map \eqref{def-AW}
into a monoidal structure on the adjoint-inverse functor $\Gam$: we define 
\begin{equation}\label{def-AW as monoidal}
 \varphi_{C,D} \ : \ \Gam C \otimes \Gam D \ \to \ \Gam(C \otimes D) 
\end{equation}
as the composite
\begin{align*}\label{monoidal structure on Gamma} 
\Gam C \otimes  \Gam D \ \xrightarrow{\eta_{\Gam C\tensor\Gam D}} \ 
\Gam N &\left(\Gam C \otimes \Gam D \right) \
\xrightarrow{\ \Gam(AW_{\Gam C,\Gam D})\ }  \\ 
&\Gam(N(\Gam C) \otimes N(\Gam D)) \ 
\xrightarrow{\Gam(\epsilon_C\tensor\epsilon_D)} \ \Gam(C \otimes D) \ . 
\end{align*}
The normalized Alexander-Whitney map \eqref{def-restricted shuffle}
is surjective (it is split by the normalized shuffle map);
since the unit and counit of the $(N,\Gam)$-adjunction
are isomorphisms, the monoidal map \eqref{def-AW as monoidal}
is also a split surjection.
The functor $\Gam$ is not lax {\em symmetric} monoidal
because the Alexander-Whitney map is not symmetric.  

So $\Gam$ induces a functor $\Gam \co \dgr \to \sring$
on the associated categories of monoids: given a differential graded ring  $R$ 
with product $\mu_R$ then $\Gam R$ is a simplicial ring with product 
\begin{equation}\label{product on Gam R} 
\Gam R \otimes \Gam R \ \xrightarrow{\varphi_{R,R}} \ \Gam (R \otimes R) 
\ \xrightarrow{\Gam(\mu_R)} \ \Gam R \ . \end{equation}
If we expand all the definitions, then the multiplication in $\Gam R$
comes out as follows: the product of two $n$-simplices
$x,y\co N\Delta^n\to R$ of $\Gam R$ is the composition
\begin{align*} N\Delta^n \ \xrightarrow{N(\text{diag})} \
N(\Delta^n\times &\Delta^n) \ \xrightarrow{AW_{\Delta^n,\Delta^n}} \\ 
& N\Delta^n\tensor N\Delta^n \
\xrightarrow{x\tensor y} \ R\tensor R \ \xrightarrow{\mu_R} \ R \ . 
\end{align*}

\begin{example}\label{ex-effect of Gam}
To give an idea of what the multiplication in $\Gam R$ looks like,
we calculate an explicit formula in the lowest dimension 
where something happens. The normalized chain complex of the
simplicial 1-simplex $\Delta^1$ has as basis the cosets of the
non-degenerate 1-simplex $\iota\in \Delta^1_1$ 
and the two vertices $0=d_1\iota$ and $1=d_0\iota$; the differential in
$N\Delta^1$ is determined by $d[\iota]=[1]-[0]$.

So for every 1-chain $r\in R_1$ of a differential graded ring $R$
we can define a chain map \mbox{$\kappa r\co N\Delta^1\to R$} by setting
\[ (\kappa r) [\iota] =r \ , \ (\kappa r)[0]=0 \quad \text{and} \quad 
(\kappa r)[1] = dr \ . \] 
This defines a monomorphism $\kappa\co R_1\to \ch(N\Delta^1,R)=(\Gam R)_1$.

The composite map
\[ N\Delta^1 \ \xrightarrow{N(\text{diag})} \
N(\Delta^1\times\Delta^1) \ \xrightarrow{AW_{\Delta^1,\Delta^1}} \
 N\Delta^1\tensor N\Delta^1 \]
is given by 
\begin{align*}
[\iota] \ &\longmapsto \ [0]\tensor\iota\, \oplus\, \iota\tensor [1] \\
[0] \ &\longmapsto \ [0]\tensor[0] \text{\quad and\quad}
[1] \ \longmapsto \ [1]\tensor[1] \ . 
\end{align*}
Hence we have
\[ (\kappa r\cdot \kappa s)[\iota] \ = \ (\kappa r)[0]\cdot (\kappa s)[\iota] 
\, + \, (\kappa r)[\iota]\cdot (\kappa s)[1] \ = \ r\cdot ds \ , \] 
and similarly $(\kappa r\cdot \kappa s)[0]=0$ and 
$(\kappa r\cdot \kappa s)[1]=dr\cdot ds$.
In other words, we have shown the formula
\[ \kappa r\cdot \kappa s\ = \ \kappa(r\cdot ds) \]
as 1-simplices of $\Gam R$, for every pair of 1-chains $r,s\in R_1$.
This formula already indicates that a simplicial ring of the
form $\Gam R$ is usually not commutative, even if $R$ is commutative in the
graded sense.
\end{example}

\subsection{Monoidal properties of $\epsilon \co N\Gamma \to \Id$
and $\eta \co \Id \to N\Gamma$}

The normalization functor is lax symmetric monoidal with structure
map induced by the shuffle map~\eqref{def-restricted shuffle}.  Thus, it also
induces a functor on the categories of monoids $N \co \sring \to \dgr$. 
We shall see in the next proposition that $N$ is left inverse to
$\Gam$ on the level of rings; however, $N$ is {\em not} right inverse
to $\Gam$ on the point-set level, but only on the
level of homotopy categories, see Remark \ref{eta not monoidal}.

\begin{lemma}\label{lem-counit is monoidal}
The adjunction counit $\epsilon\co N\Gam\to \text{\em Id}_{\ch}$
is a monoidal transformation with respect to the composite
monoidal structure on $N\Gam$. More precisely,
for every pair of connective chain complexes $C$ and $D$,
the following diagram commutes
\begin{equation}\label{eq-counit monoidal} 
\xymatrix@C=20mm{ N\Gam C\tensor N\Gam D \ar^{\nabla_{\Gam C,\Gam D}}[r] 
\ar_-{\epsilon_C\tensor\epsilon_D}[drr] &
N(\Gam C \otimes \Gam D) \ar^{N(\varphi_{C,D})}[r] & 
N \Gam(C \otimes D) \ar^{\epsilon_{C\tensor D}}[d] \\
& & C\otimes D }\end{equation}
\end{lemma}
\begin{proof}
The proof is a diagram chase, the main ingredient of which is the fact 
that the composite \eqref{AW nabla=Id} of
the normalized shuffle and Alexander-Whitney maps is the identity.
We start with the identity
\[ \nabla_{\Gam C,\Gam D}\, \circ\, \epsilon_{N\Gam C\tensor N\Gam D} \ = \ 
\epsilon_{N(\Gam C\tensor \Gam D)}\,\circ\, N\Gam(\nabla_{\Gam C,\Gam D}) \]
as morphisms from $N\Gam(N\Gam C\tensor N\Gam D)$ to 
$N(\Gam C\tensor \Gam D)$, which just says that $\epsilon$ is natural.
The map $\epsilon_{N(\Gam C\tensor \Gam D)}$ is inverse
to the map $N(\eta_{\Gam C\tensor\Gam D})\co 
N(\Gam C \otimes \Gam D)\to N\Gam N(\Gam C \otimes \Gam D)$,
so we can rewrite the previous identity as
\[ N(\eta_{\Gam C\tensor\Gam D})\, \circ\, 
\nabla_{\Gam C,\Gam D}\,\circ\, \epsilon_{N\Gam C\tensor N\Gam D} \ = \ 
N\Gam(\nabla_{\Gam C,\Gam D}) \]
as morphisms from $N\Gam(N\Gam C\tensor N\Gam D)$ to 
$N\Gam N(\Gam C\tensor \Gam D)$.
Now we compose with the map $N\Gam(AW_{\Gam C,\Gam D})\co 
N\Gam N(\Gam C\tensor \Gam D)\to N\Gam (N\Gam C\tensor N\Gam D)$
and exploit that the Alexander-Whitney map is left inverse to the
shuffle map (see \eqref{AW nabla=Id}); this yields
\begin{align*}  N\Gam(AW_{\Gam C,\Gam D}&)\,\circ\, 
N(\eta_{\Gam C\tensor\Gam D})\,\circ\, 
\nabla_{\Gam C,\Gam D}\,\circ\, \epsilon_{N\Gam C\tensor N\Gam D} \\
&= \ N\Gam(AW_{\Gam C,\Gam D})\,\circ\, N\Gam(\nabla_{\Gam C,\Gam D}) \
= \ \Id_{N\Gam (N\Gam C\tensor N\Gam D)} \ .   \end{align*} 
Composing with $N\Gam(\epsilon_C\tensor\epsilon_D)\co 
N\Gam (N\Gam C\tensor N\Gam D)\to N\Gam(C\tensor D)$ and
substituting the definition \eqref{def-AW as monoidal}
of the monoidal transformation $\varphi$, we get
\begin{align*}
N\Gam(\epsilon_C\tensor\epsilon_D) \ &= \ \\
N\Gam(\epsilon_C&\tensor\epsilon_D)
\,\circ\,  N\Gam(AW_{\Gam C,\Gam D})\,\circ\, 
N(\eta_{\Gam C\tensor\Gam D})\,\circ\, 
\nabla_{\Gam C,\Gam D}\,\circ\, \epsilon_{N\Gam C\tensor N\Gam D} \\
& = \ N(\varphi_{C,D})\, \circ\, 
\nabla_{\Gam C,\Gam D}\,\circ\, \epsilon_{N\Gam C\tensor N\Gam D} \ . 
\end{align*} 
Since the counit $\epsilon$ is invertible, we can rewrite this as
\begin{eqnarray*}\epsilon_{C\tensor D}\circ N(\varphi_{C,D}) 
\circ \nabla_{\Gam C,\Gam D}
& = & \epsilon_{C\tensor D}\circ N\Gam(\epsilon_C\tensor\epsilon_D)
\circ \epsilon^{-1}_{N\Gam C\tensor N\Gam D} 
\ = \ \epsilon_C\tensor \epsilon_D 
\end{eqnarray*} 
(the second equation is the naturality of $\epsilon$). 
This final identity is saying that the transformation $\epsilon\co N\Gam\to\Id$
is monoidal.
\end{proof}

\begin{proposition}\label{prop-full embedding} The functor
\[ \Gam \ : \ \dgr \ \to \ \sring \] 
which sends a connective differential graded ring $R$ to the
simplicial abelian group $\Gam R$ with multiplication \eqref{product on Gam R}
is full and faithful.
The composite endo-functor $N\Gam$ of the category of
differential graded rings is naturally isomorphic to the identity functor.
\end{proposition}
\begin{proof}
The algebra valued functor $\Gam$ is induced from an equivalence 
between the underlying categories of simplicial abelian groups
and chain complexes.
So in order to show that $\Gam$ is fully faithful we have to prove that
for every morphism $f\co \Gam R\to \Gam S$ of simplicial rings,
the unique morphism $g\co R\to S$ of chain complexes
which satisfies $f=\Gam(g)$ is multiplicative and preserves the units.

The non-trivial part is to show that if $f\co \Gam R\to \Gam S$ 
is multiplicative, then the unique preimage $g\co R\to S$
is also  multiplicative. Since $\Gam$ is an equivalence 
on underlying categories, in order to show the relation
$\mu_S\circ(g\tensor g)=g\circ \mu_R$ as  chain maps from $R\tensor R$ to $S$, 
we may as well show the relation
$\Gam(\mu_S\circ(g\tensor g))=\Gam(g\circ \mu_R)$
as maps of simplicial abelian groups from $\Gam(R\tensor R)$ to $\Gam S$.
In the diagram
\[\xymatrix@C=15mm{ \Gam R \tensor \Gam R \ar[r]^-{\varphi_{R,R}} 
\ar[d]_{\Gam(g)\tensor\Gam(g)} & 
\Gam(R\tensor R) \ar[r]^-{\Gam(\mu_R)} \ar[d]_{\Gam(g\tensor g)} &
\Gam R \ar[d]^{\Gam(g)} \\
\Gam S\tensor \Gam S \ar[r]_-{\varphi_{S,S}} & 
\Gam(S\tensor S) \ar[r]_-{\Gam(\mu_S)} & \Gam S }\]
the left square commutes since $\varphi$ is natural, and the composite
square commutes since $\Gam(g)=f$ is multiplicative. Since the
upper left morphism $\varphi_{R,R}$ is surjective, 
the right square commutes as well.

The natural isomorphism $N\Gam\iso\Id$
is given by the counit $\epsilon\co N\Gam\to \Id$
of the adjunction-equivalence between $N$ and $\Gam$.
This counit is an isomorphism, and it is monoidal by 
Lemma \ref{lem-counit is monoidal}; this implies that on ring objects,
the map $\epsilon$ is multiplicative; that $\epsilon$ is unital is even easier,
so $\epsilon$ is a natural isomorphism of connective differential graded rings
when evaluated on such objects.
\end{proof}

\begin{remark}\label{eta not monoidal}
The unit $\eta\co \Id \to \Gam N$ of the adjunction-equivalence 
between $N$ and $\Gam$ is {\em not} monoidal.
More precisely, the composite
\begin{align}\label{eq-eta composite}
A\tensor B \ \xrightarrow{\eta_A\tensor\eta_B} \
\Gam (NA) \otimes \Gam &(NB)  \ \xrightarrow{\varphi_{NA,NB}} \\
&\Gam(NA \otimes NB) \ \xrightarrow{\Gam(\nabla_{A,B})}\ \Gam N(A\otimes B)
\nonumber\end{align}
need not in general be equal to the map 
$\eta_{A\tensor B}\co  A\tensor B\to\Gam N(A\otimes B)$.
(Consider for example $A = B = 
\widetilde{\mathbb{Z}} (\Delta[1]/\partial\Delta[1])=\Gam(\mathbb{Z}[1])$.  
In dimension one the composite \eqref{eq-eta composite} is zero 
since it factors through $\Gam(\mathbb{Z}[2])$. But in dimension one, 
$\eta_{A\tensor B}$ is an isomorphism between free abelian groups of rank two.)
The situation is worse than for the counit $\epsilon$ 
(compare Lemma \ref{lem-counit is monoidal}) 
because the composite of normalized shuffle and Alexander-Whitney map
in the other order is only homotopic, but not equal to, the identity.
Correspondingly, the composite \eqref{eq-eta composite} is homotopic, 
but not necessarily equal to, the map $\eta_{A\tensor B}$.

Nevertheless, the composite $\Gam N$ is connected by a chain of
two natural weak equivalences to the identity functor on the category of
simplicial rings. In order to see this though, we have to refer
to the Quillen equivalence of Theorem~\ref{thm-sab-ch}~\eqref{iii}.
When considered as a ring valued functor, $N\co \sring\to\dgr$
has a left adjoint $L^{\text{mon}}\co \dgr\to\sring$,
see Section \ref{sub-left adjoint}, which is {\em not} given by $\Gam$ 
on underlying chain complexes. 
Moreover, the adjoint pair $N$ and $L^{\text{mon}}$
form a Quillen equivalence (Theorem \ref{thm-sab-ch} \eqref{iii} for $k=\mZ$).

Given a simplicial ring $A$, we choose a cofibrant replacement
\[ q \ : \ (NA)^c \ \xrightarrow{\sim} \ NA \] 
of $NA$ in the model category of connective differential graded rings
(see Section \ref{sub-model structures}). 
The model structure of differential graded
rings is cofibrantly generated, so the small object argument provides
a {\em functorial} choice of such a cofibrant replacement. 

Since $N$ and $L^{\text{mon}}$ are a Quillen equivalence, the adjoint morphism 
\[ \tilde q\ : \ L^{\text{mon}}((NA)^c)\ \to\ A \] 
is a weak equivalence of simplicial rings. 
By Lemma \ref{lem-counit is monoidal}, the adjunction counit
$\epsilon_{NA}\co N\Gam NA\to NA$ is an isomorphism 
of differential graded rings.
So we can form the composite multiplicative quasi-isomorphism 
$\epsilon_{NA}^{-1}\circ q\co (NA)^c\to N\Gam NA$
and take its monoid-valued adjoint
\begin{equation}\label{eq-monoid adjoint} 
L^{\text{mon}}((NA)^c) \ \to \ \Gam NA \ . \end{equation}
Since we have a Quillen equivalence, \eqref{eq-monoid adjoint} is also 
a weak equivalence of simplicial rings. 
Altogether we obtain a chain of natural weak equivalences of simplicial rings
\[ \Gam NA \ \xleftarrow{\ \sim\ } \ L^{\text{mon}}((NA)^c) \ 
\xrightarrow{\ \sim\ } \ A \ . \]
It is tempting to add the adjunction unit 
$\eta_A\co \Gam NA\to A$ to directly connect the two simplicial rings; 
but $\eta_A$ is not in general a multiplicative map, and the resulting triangle
involving $\Gam NA, L^{\text{mon}}((NA)^c)$ and $A$ need not commute !
\end{remark}

\section{Weak monoidal equivalences}\label{sec-main}

In this section we first discuss the definitions of monoidal
structures and their interactions with model category structures.
Section~\ref{3.1} recalls the notion of a {\em monoidal
model category} which is a model category with a compatible monoidal
product.  In Section~\ref{3.2} we define the notions of {\em weak} and 
{\em strong  monoidal Quillen equivalences} between two monoidal 
model categories.  A weak monoidal Quillen equivalence provides the basic 
properties necessary for lifting the Quillen equivalence to
categories of monoids and modules.  In a weak monoidal Quillen equivalence
the right adjoint is assumed to be lax monoidal and hence induces
functors on the associated categories of monoids and modules.  This
is not assumed for the left adjoint, though.  Section~\ref{sub-left adjoint}
discusses the induced right adjoints 
and the relationship between the various context-dependent left adjoints.
With this background we can then state our main result, 
Theorem~\ref{thm-WMQE gives equivalent algebras},
about Quillen equivalences on categories of monoids and modules.
In Section~\ref{subsec-criterion} we discuss a criterion for establishing when
a Quillen adjoint pair is a weak monoidal Quillen pair.  The general
criterion is given in Proposition~\ref{prop-generating criterion}; 
a variant for stable model categories appears in Proposition~\ref{prop-stable
generating criterion}. 

\subsection{Monoidal model categories}\label{3.1}
We consider a closed symmetric monoidal category~\cite[6.1]{bor} $\cC$
and we denote the monoidal product by $\otimes$ (sometimes by $\sm$), 
the unit object by $\mIC$ and the internal function objects by 
$\Hom_{\cC}(-,-)$.
The internal function objects are almost never used 
explicitly (but see the proof of Proposition \ref{prop-generating criterion}).
But having a right adjoint makes sure 
that the monoidal product preserves colimits in both variables.

\begin{definition}\label{def-mmc}
A model category $\cC$ is a {\em monoidal model category} if it has 
a closed symmetric monoidal structure with product $\tensor$ and
unit object $\mI$ and satisfies the following two axioms.

{\bf Pushout product axiom}\qua 
Let $A \xrightarrow{} B$ and $K \xrightarrow{} L$ be cofibrations in $\cC$.
Then the map
\[ A\tensor L \, \amalg_{A\tensor K} \, B\tensor K \ \to \ B\tensor L \]
is also a cofibration.  If in addition one of the former maps is a weak
equivalence, so is the latter map.

{\bf Unit axiom}\qua Let $q\co \mI^c\stackrel{\sim}{\to} \mI$ 
be a cofibrant replacement of the unit object. 
Then for every cofibrant object $A$, the morphism
$q\tensor\Id\co \mI^c\tensor A\to \mI\tensor A\iso A$ is a weak equivalence.
\end{definition}

The previous definition is essentially the same as that of a 
{\em symmetric monoidal model category} in~\cite[4.2.6]{hovey-book};
the only difference is that a model category in Hovey's sense is also equipped
with a choice of cofibrant replacement functor, and Hovey requires 
the unit axiom for the particular functorial cofibrant replacement 
of the unit object. But given the pushout product axiom, then the
unit axiom holds for one choice of cofibrant replacement if and only
if it holds for any other choice. Of course the unit axiom is redundant
if the unit object is cofibrant, as we often assume in this paper.
The unit axiom did not occur in the definition of a monoidal model
category in~\cite{ss}, since it did not play a role in the arguments
of that paper.

In this paper we are interested in model categories 
of monoids~\cite[VII 3]{ML-working} and modules
(i.e., objects with an action by a monoid, see~\cite[VII 4]{ML-working})
in some underlying monoidal model category. 
In the cases we study the model structure is always transferred 
or lifted from the underlying category to the category 
of more structured objects as in the following definition.

\begin{definition}\label{def-create} 
Consider a functor $R\co \cT\to \cC$ to a model category $\cC$ 
with a left adjoint $L\co \cC\to \cT$. We call an object of $\cT$ a 
{\em cell object} if it can be obtained from the initial object 
as a (possibly transfinite) composition~\cite[2.1.1]{hovey-book}
of pushouts along morphisms of the form $Lf$, for $f$ a cofibration in $\cC$.
We say that the functor $R$ {\em creates a model structure} if 
\begin{itemize}
\item the category $\cT$ supports a model structure 
(necessarily unique) in which a morphism $f\co X\to Y$ 
is a weak equivalence, respectively fibration,
if and only the morphism $Rf$ is a weak equivalence, 
respectively fibration, in $\cC$, and
\item every cofibrant object in $\cT$ is a retract of a cell object.
\end{itemize}
\end{definition}

The typical example of Definition \ref{def-create} occurs when
the model structure on $\cC$ is cofibrantly generated  
(\cite[2.2]{ss} or~\cite[2.1.17]{hovey-book}) and then a lifting theorem 
is used to lift the model structure to $\cT$ along the adjoint functor pair.
In~\cite[4.1]{ss} we give sufficient conditions for ensuring 
that the forgetful functors to an underlying monoidal model category $\cC$ 
create model structures for monoids or modules and algebras:
it suffices  that the model structure on $\cC$ 
is cofibrantly generated, that the objects of $\cC$ are small 
relative to the whole category (\cite[\S 2]{ss} or~\cite[2.1.3]{hovey-book}) 
and that the {\em monoid axiom}~\cite[3.3]{ss} holds.

Another situation where the forgetful functor creates a model structure
for the category of monoids is when all objects in $\cC$ are fibrant and 
when there exists an `interval with coassociative, comultiplication';
for more details compare the example 
involving chain complexes~\cite[2.3 (2), Sec.\ 5]{ss}
or more generally~\cite[Prop. 4.1]{Berger-Moerdijk}.

\subsection{Monoidal Quillen pairs}\label{3.2}

The main goal for this paper is to give conditions which show that a Quillen 
equivalence between two monoidal model categories 
induces a Quillen equivalence 
on the categories of monoids, modules and algebras. 
For this we have to assume that the functors involved preserve
the monoidal structure in some way.  We assume that the right adjoint
is lax monoidal in the sense of the following definition.

\begin{definition}\label{lax}
A {\em lax monoidal} functor between 
monoidal categories is a functor $R\co \cC\to \cD$ equipped with 
a morphism $\nu\co \mID \to R(\mIC)$ and natural morphisms
\[  \varphi_{X,Y} \ \co \ RX \sm RY \ \to \ R(X\otimes Y)  \]
which are coherently associative and unital
(see diagrams 6.27 and 6.28 of \cite{bor}).
A lax monoidal functor is {\em strong monoidal} if the morphisms
$\nu$ and $\varphi_{X,Y}$ are isomorphisms.
\end{definition}

Consider a lax monoidal functor $R\co \cC\to \cD$ between monoidal categories,
with monoidal structure maps $\nu$ and $\varphi_{X,Y}$.
If $R$ has a left adjoint $\lambda\co \cD\to\cC$, we can consider 
the adjoint $\tilde\nu\co \lambda(\mID)\to\mIC$ of $\nu$ and the natural map 
\begin{equation}\label{eq-tvarphi}
\tvarphi\co\lambda( A \sm B) \ \to \ \lambda A \otimes \lambda B \end{equation}
adjoint to the composite
\[ A \sm B \ \xrightarrow{\eta_A\sm\eta_B} \  
R\lambda A \sm R \lambda B \ \xrightarrow{\varphi_{\lambda A,\lambda B}} 
\ R(\lambda A \otimes \lambda B) \ . \]
The map $\tvarphi$ can equivalently be defined as the composition
\begin{align}\label{eq-second def of tvarphi} 
\lambda( A \sm B) \ \xrightarrow{\lambda(\eta_A\sm\eta_B)} \ 
\lambda( R\lambda A &\sm R\lambda B) 
\xrightarrow{\lambda(\varphi_{\lambda_A,\lambda_B})} \\ 
&\lambda R(\lambda A \otimes \lambda B) 
\ \xrightarrow{\epsilon_{\lambda_A\tensor\lambda_B}} \
\lambda A \otimes \lambda B \ ; \nonumber\end{align} 
here $\eta$ and $\epsilon$ denote the unit respectively counit
of the adjunction.
With respect to these maps, $\lambda$ is a {\em lax comonoidal} functor
(also referred to as an {\em op-lax monoidal} functor). 
The map $\tvarphi$ need not be an isomorphism; 
in that case $\lambda$ does not have a {\em monoidal} structure, 
and so it does not pass to a functor on the monoid and module categories.

\begin{definition}\label{def-WMQ equivalence} 
A pair of adjoint functors 
\[\xymatrix@=12mm{ \cD \quad \ar@<-.4ex>_-{\lambda}[r] &
\quad  \cC \quad  \ar@<-.4ex>_-{R}[l] } 
\]
between model categories is a {\em Quillen adjoint pair} 
if the right adjoint $R$ preserves fibrations and trivial fibrations.  
A Quillen adjoint pair induces adjoint total derived functors
between the homotopy categories by~\cite[I.4.5]{Q}.
A Quillen functor pair is a {\em Quillen equivalence} if the total
derived functors are adjoint equivalences of the homotopy categories.

A {\em weak monoidal Quillen pair} between monoidal model categories
$\cC$ and $\cD$ consists of a Quillen adjoint functor pair 
$(\!\xymatrix{\lambda\co \cD \ar@<-.4ex>[r] & \cC : R\ar@<-.4ex>[l]}\!)$
with a lax monoidal structure on the right adjoint
\[  \varphi_{X,Y} \ : \ RX \sm RY \ \to \ R(X\otimes Y)  \ , \quad 
\nu \ : \ \mID \ \to \ R(\mIC) \]
such that the following two conditions hold:
\begin{itemize}
\item[(i)] for all cofibrant objects $A$ and $B$ in $\cD$ the 
comonoidal map \eqref{eq-tvarphi}
\[ \tvarphi \ \co\ \lambda( A \sm B) \ \to \ \lambda A \otimes \lambda B \]
is a weak equivalence in $\cC$ and 
\item[(ii)] for some (hence any) cofibrant replacement
$q\co\mID^c\xrightarrow{\sim}\mID$ of the unit object in $\cD$, the composite map
\[ \lambda(\mID^c)\ \xrightarrow{\, \lambda(q)\, }\ 
\lambda(\mID)\ \xrightarrow{\ \tilde\nu\ } \ \mIC \] 
is a weak equivalence in $\cC$.
\end{itemize}

A {\em strong monoidal Quillen pair} is a weak monoidal Quillen pair
for which the comonoidal maps $\tvarphi$ and $\tilde\nu$ are isomorphisms.
Note that if $\mID$ is cofibrant and $\lambda$ is strong monoidal, then
$R$ is lax monoidal and the Quillen pair is a strong monoidal Quillen pair.

A weak (respectively strong) monoidal Quillen pair is a 
{\em weak monoidal Quillen equivalence}
(respectively {\em strong monoidal Quillen equivalence})
if the underlying Quillen pair is a Quillen equivalence.
\end{definition}

Strong monoidal Quillen pairs are the same as
{\em monoidal Quillen adjunctions} in the sense of 
Hovey~\cite[4.2.16]{hovey-book}. 
The weak monoidal Quillen pairs do not occur in Hovey's book.

\subsection{Various left adjoints}\label{sub-left adjoint} 
As any lax monoidal functor, the right adjoint $R\co \cC\to\cD$ of a weak monoidal
Quillen pair induces various functors on the categories of monoids and modules.
More precisely, for a monoid $A$ in $\cC$ 
with multiplication $\mu\co A\tensor A\to A$ and unit $\eta\co \mIC\to A$,
the monoid structure on $RA$ is given by the composite maps 
\[ RA \sm RA \ \xrightarrow{\ \varphi_{A,A}} \  R(A \otimes A) 
\ \xrightarrow{R(\mu)} \ RA \  \text{\quad and \quad} 
\mID \ \xrightarrow{\ \nu\ } \ R(\mIC)  \ \xrightarrow{R(\eta)} \ RA  \ . \] 
Similarly, for an $A$-module $M$ with action morphism
$\alpha\co M\tensor A\to M$, the $\cD$-object $RM$ becomes
an $RA$-module via the composite morphism
\[ RM \sm RA \ \xrightarrow{\ \varphi_{M,A}} \  R(M \otimes A) 
\ \xrightarrow{R(\alpha)} \ RM \ . \]

In our context, $R$ has a left adjoint $\lambda\co \cD\to\cC$.
The left adjoint inherits an `adjoint' {\em co}monoidal structure
$\tvarphi \co\lambda(A \sm B)\to\lambda A \otimes \lambda B$ and
$\tilde\nu\co \lambda(\mID)\to\mIC$, see~\eqref{eq-tvarphi}, 
and the pair is strong monoidal if $\tvarphi$ and $\tilde\nu$ are isomorphisms.
In that case, the left adjoint becomes a
strong {\em monoidal} functor via the inverses
\[ \tvarphi^{-1} \co \lambda A \otimes \lambda B \ \to \ \lambda(A \sm B)\ 
\text{\quad and \quad} 
\tilde\nu^{-1} \co \mIC\ \to \ \lambda(\mID) \ . \]
Via these maps, $\lambda$ then lifts to a functor 
on monoids and modules, and the lift is again adjoint
to the module- or algebra-valued version of $R$.

However, we want to treat the more general situation of
{\em weak} monoidal Quillen pairs.
In that case, the functors induced by $R$ on modules and algebras
still have left adjoints.
However, on underlying $\cD$-objects, these left adjoints 
are {\em not} usually given by the original left adjoint $\lambda$. 
As far as $R$ is concerned, we allow ourselves the
abuse of notation to use the same symbol for the original 
lax monoidal functor $R$ from $\cC$ to $\cD$ as well as for its structured
versions. However, for the left adjoints it seems more appropriate to
use different symbols, which we now introduce.

In our applications we always assume that the forgetful functor
from $\cC$-monoids to $\cC$ creates a model structure.
In particular, the category of $\cC$-algebras has colimits and 
the forgetful functor has a left adjoint `free monoid' (or `tensor algebra') 
functor~\cite[VII 3, Thm.~2]{ML-working}
\[ T_{\cC}X \ = \ \coprod_{n\geq 0} \, X^{\tensor n} \ = \ 
\mIC \, \amalg \, X \, \amalg \, (X\tensor X) \, \amalg \,  \dots \ , \]
with multiplication given by juxtaposition, and similarly for $T_{\cD}$.
This implies that the monoid-valued lift $R\co \cC\mon \to \cD\mon$
again has a left adjoint
\[ L^{\text{mon}} \ : \ \cD\mon \ \to \ \cC\mon \ . \] 
Indeed, for a $\cD$-monoid $B$, the value of the left adjoint
can be defined as the coequalizer of the two $\cC$-monoid morphisms
\[\xymatrix@=12mm{ T_{\cC}\left( \lambda(T_{\cD}B)\right) 
\ \ar@<-.4ex>[r] \ar@<.4ex>[r] & \ T_{\cC}( \lambda B) \ \ar[r] &
L^{\text{mon}}B  } 
\]
(where the forgetful functors are not displayed).
One of the two maps is obtained from the adjunction unit $T_{\cD}B\to B$
by applying the the composite functor $T_{\cC}\lambda$; the other
map is the unique $\cC$-monoid morphism which restricts to the $\cC$-morphism
\[ \lambda(T_{\cD}B) \ \iso \ \coprod_{n\geq 0} \, 
\lambda\left(B^{\sm n}\right) \ \xrightarrow{\ \amalg\tvarphi \ } \  
\coprod_{n\geq 0} \, (\lambda B)^{\tensor n} \ \iso \ 
T_{\cC}( \lambda B) \ . \] 
Since $R$ preserves the underlying objects, the monoid left adjoint
and the original left adjoint are related 
via a natural isomorphism of functors from $\cD$ to monoids in $\cC$
\begin{equation}\label{eq-free monoid iso}
L^{\text{mon}}\circ T_{\cD} \ \iso \ T_{\cC}\circ \lambda \ .
\end{equation}

As in the above case of monoids, the module valued functor 
$R\co \Modr A\to\Modr RA$ has a left adjoint
\[ L^A \ : \ \Modr RA \ \to \ \Modr A \] 
as soon as free $R$-modules and coequalizers of $R$-modules exist.
Since $R$ preserves the underlying objects, this module left adjoint
and the original left adjoint are related 
via a natural isomorphism of functors from $\cD$ to $A$-modules
\begin{equation}\label{eq-free module iso for L^A}
L^A\circ(-\sm RA) \ \iso \ (-\tensor A)\circ \lambda  \ ; \end{equation}
here $X\tensor A$ is the free $A$-module generated by a $\cC$-object $X$,
and similarly $Y\sm RA$ for a $\cD$-object $Y$.

Finally, for a monoid $B$ in $\cD$, the lax monoidal functor $R$
induces a functor from the category of $L^{\text{mon}}B$-modules
to the category of $B$-modules; this is really the composite functor
\[ \Modr (L^{\text{mon}}B) \ \xrightarrow{\ R\ }  \ \Modr R(L^{\text{mon}}B)
\ \xrightarrow{\ \eta^*} \Modr B \]
where the second functor is restriction of scalars along the
monoid homomorphism (the adjunction unit) $\eta \co  B\to R(L^{\text{mon}}B)$.
We denote by
\begin{equation}\label{eq-L_B} 
L_B \ : \ \Modr B \ \to \ \Modr (L^{\text{mon}}B) \end{equation}
the left adjoint to the functor $R\co \Modr (L^{\text{mon}}B)\to\Modr B$.
This left adjoint factors as a composition
\[ \Modr B \ \xrightarrow{\ -\tensor_{B}R(L^{\text{mon}}B)} \
\Modr R(L^{\text{mon}}B) \ \xrightarrow{\ L^{L^{\text{mon}}B}} \ 
\Modr (L^{\text{mon}}B)  \]
(the first functor is extension of scalars along $B\to R(L^{\text{mon}}B)$);
the left adjoint is related to the free module functors by a
natural isomorphism
\begin{equation}\label{eq-free module iso for L_B}
L_B\circ(-\sm B) \ \iso \ (-\tensor L^{\text{mon}}B)\circ \lambda
\ . \end{equation}

We repeat that if the monoidal pair is {\em strong monoidal},
then left adjoints $L^{\text{mon}}$ and $L_B$ are
given by the original left adjoint $\lambda$, which is
then monoidal via the inverse of $\tvarphi$.
Moreover, the left adjoint $L^A\co \Modr RA\to\Modr A$ is then given 
by the formula
\[ L^A(M) \ = \ \lambda(M)\sm_{\lambda(RA)} A \ , \] 
where $A$ is a $\lambda (RA)$-module via the adjunction counit
$\epsilon\co \lambda (RA)\to A$.
In general however, $\lambda$ does not pass to monoids and modules, 
and the difference between $\lambda$ and the structured
adjoints $L^{\text{mon}}$, $L^A$ and $L_B$ is investigated
in Proposition \ref{prop-L} below.

We need just one more definition before stating our main theorem.

\begin{definition}\label{def-Quillen invariance}
Let $(\cC, \tensor, \mIC)$ be a monoidal model category such that the
forgetful functors create model structures for modules over any monoid. 
We say that {\em Quillen invariance holds}
for $\cC$ if for every 
weak equivalence of $\cC$-monoids $f\co R\to S$,
restriction and extension of scalars along $f$ induce a Quillen equivalence
between the respective module categories.
\[\xymatrix{- \sm_{R} S \ : \ \Modr (R)\ \ar@<-.4ex>[r]&
\ \Modr S \ : f^* \ar@<-.4ex>[l] }\]
\end{definition}

A sufficient condition for Quillen invariance in $\cC$ is that for 
every cofibrant right $R$-module $M$ the functor $M\tensor_R-$ 
takes weak equivalences of left $R$-modules to weak equivalences in $\cC$
(see for example~\cite[4.3]{ss} or Theorem~\ref{a.1.2} (2)).

\begin{theorem}\label{thm-WMQE gives equivalent algebras}
Let $R\co \cC\to\cD$ be the right adjoint of a 
weak monoidal Quillen equivalence. Suppose that the unit objects in
$\cC$ and $\cD$ are cofibrant.

\begin{enumerate}
\item\label{ONE} Consider a cofibrant monoid $B$ in $\cD$ such that
the forgetful functors create model structures
for modules over $B$ and modules over $L^{\text{\em mon}}B$.
Then the adjoint functor pair 
\[\xymatrix{L_B\ : \ \Modr B \ar@<-.4ex>[r]&
\ \Modr (L^{\text{\em mon}}B)\ : \ R\ar@<-.4ex>[l] }\]
is a Quillen equivalence.

\item\label{TWO}
Suppose that Quillen invariance holds in $\cC$ and $\cD$. 
Then for any fibrant monoid $A$ in $\cC$ such that the forgetful functors 
create model structures for modules over $A$ and modules over $RA$,
the adjoint functor pair
\[\xymatrix{L^A\ : \ \Modr RA\ \ar@<-.4ex>[r]&\ \Modr A \ :\ 
R\ar@<-.4ex>[l] }\]
is a Quillen equivalence.
If the right adjoint $R$ preserves weak equivalences between
monoids and the forgetful functors create model structures
for modules over any monoid, then this holds for any monoid $A$ in $\cC$.

\item\label{THREE} If the forgetful functors create model structures
for monoids in $\cC$ and $\cD$, then the adjoint functor pair
\[\xymatrix{L^{\text{\em mon}}\ : \ \cD\mon\ \ar@<-.4ex>[r]&
\ \cC\mon \ : \ R\ar@<-.4ex>[l] }\]
is a Quillen equivalence between the model categories of monoids.
\end{enumerate}
\end{theorem}

The statements \eqref{ONE} and \eqref{TWO} for modules 
in the previous Theorem \ref{thm-WMQE gives equivalent algebras} generalize 
to `rings with many objects' or enriched categories, 
see Theorem~\ref{thm-WMQE many generators}.
The proof of Theorem~\ref{thm-WMQE gives equivalent algebras} 
appears in Section~\ref{sec-proofs}. 

\subsection{A Criterion for weak monoidal pairs}
\label{subsec-criterion}

In this section we assume that $\lambda$ and $R$ 
form a Quillen adjoint functor pair
between two monoidal model categories 
$(\cC, \otimes, \mIC)$ and  $(\cD, \sm, \mID )$.
We establish a sufficient condition for when the
Quillen pair is a weak monoidal Quillen pair:
if the unit object $\mID$ detects weak equivalences
(see Definition~\ref{def-detect}), then
the lax comonoidal transformation \eqref{eq-tvarphi}
\[ \tvarphi\co\lambda( A \sm B) \ \to \ \lambda A \otimes \lambda B \ . \]
is a weak equivalence on cofibrant objects.

To define what it means to detect weak equivalences,
we use the notion of a {\em cosimplicial resolution}
which was introduced by Dwyer and Kan~\cite[4.3]{DK} as a device 
to provide homotopy meaningful mapping spaces. 
More recently, \mbox{(co-)}simplicial resolutions have been called
{\em \mbox{(co-)}simplicial frames}~\cite[5.2.7]{hovey-book}.
Cosimplicial objects in any model category admit the    
{\em Reedy model structure} in which the weak equivalences are 
the cosimplicial maps which are levelwise weak equivalences 
and the cofibrations are the {\em Reedy cofibrations}; 
the latter are a special class of levelwise cofibrations defined with
the use of `latching objects'. 
The Reedy fibrations are defined by the
right lifting property for Reedy trivial cofibrations or equivalently
with the use of {\em matching objects}; see \cite[5.2.5]{hovey-book}
for details on the Reedy model structure.
A {\em cosimplicial resolution} of an object $A$ of $\cC$ is a cofibrant
replacement $A^*\to cA$ in the Reedy model structure
of the constant cosimplicial object $cA$ with value $A$.
In other words, a cosimplicial resolution is a Reedy
cofibrant cosimplicial object which is {\em homotopically constant} 
in the sense that each cosimplicial structure map 
is a weak equivalence in $\cC$.
Cosimplicial resolutions always exist~\cite[4.5]{DK} and are unique 
up to level equivalence under $A$.

If $A^*$ is a cosimplicial object and $Y$ is an object of $\cC$, then there
is a simplicial set $\map(A^*,Y)$ of $\cC$-morphisms defined by
\[ \map(A^*,Y)_n \ = \ {\cC}(A^n,Y) \ . \]
If $f\co A^*\to B^*$ is a level equivalence between Reedy 
cofibrant cosimplicial 
objects and $Y$ is a fibrant object, then the induced map of mapping spaces 
$\map(f,Y)\co \map(B^*,Y)\to \map(A^*,Y)$ 
is a weak equivalence~\cite[5.4.8]{hovey-book}. 
Hence the homotopy type of the simplicial set $\map(A^*,Y)$,
for $A^*$ a cosimplicial resolution, depends only on the underlying object
$A^0$ of $\cC$. The path components of the simplicial set $\map(A^*,Y)$
are in natural bijection with the set of homotopy classes of maps 
from $A^0$ to $Y$~\cite[5.4.9]{hovey-book},
\[ \pi_0 \map(A^*,Y) \ \iso \ [A^0,Y]_{\Ho(\cC)} \ . \]

\begin{remark} The notion of a cosimplicial resolution is modeled on
the `product with simplices'. More precisely, 
in a {\em simplicial} model category $\cC$, we have 
a pairing between objects of $\cC$ and simplicial sets.
So if we let $n$ vary in the simplicial category $\Delta$, 
we get a cosimplicial object $\Delta^n\times A$ associated to
every object $A$ of $\cC$. If $A$ is cofibrant, then this cosimplicial
object is Reedy cofibrant and homotopically constant,
i.e., a functorial cosimplicial resolution of $A$.
Moreover, the $n$-simplices of the simplicial set $\map(A,Y)$ 
--- which is also part of the simplicial structure --- are in natural 
bijection with the set of morphisms from $\Delta^n\times A$ to $Y$,
\[ \map(A,Y)_n \ = \ {\cC}(\Delta^n\times A,Y) \ . \] 
So the mapping spaces with respect to the preferred 
resolution $\Delta^\ast\times A$ coincide 
with the usual simplicial function spaces~\cite[II.1.3]{Q}.
\end{remark}

\begin{definition}\label{def-detect}
An object $A$ of a model category $\cC$ {\em detects weak equivalences} if
 for some (hence any) cosimplicial frame $A^*$ of $A$ 
the following condition holds:
a morphism $f\co Y \to Z$ between fibrant objects is a weak equivalence 
if and only if the map
\[ \map(A^*,f) \ : \ \map(A^*,Y)\ \to\  \map(A^*,Z) \]
is a weak equivalence of simplicial sets.
\end{definition}

\begin{example}\label{ex-free rank one generates}
 Every one point space detects weak equivalences
in the model category of topological spaces with respect 
to weak homotopy equivalences. The one point simplicial set 
detects weak equivalences of simplicial sets with respect 
to weak homotopy equivalences.
For a simplicial ring $R$, the free $R$-module of rank one detects 
weak equivalences of simplicial $R$-modules.

Let $A$ be a connective differential graded algebra. Then the free 
differential graded $A$-module of rank one detects weak equivalences of
connective $A$-modules. Indeed, a cosimplicial resolution of the  
free $A$-module of rank one is given by $N\Delta^*\tensor A$,
i.e., by tensoring $A$ with the normalized chain complexes of the
standard simplices. With respect to this resolution the mapping space
into another connective $A$-module $Y$ has the form
\[ \map_{\Modr A}(N\Delta^*\tensor A,Y) \ \iso \
 \map_{\ch}(N\Delta^*,Y) \ = \ \Gamma Y \ . \] 
A homomorphism of connective $A$-modules is a quasi-isomorphism
if and only if it becomes a weak equivalences of simplicial sets
after applying the functor $\Gam$. 
Hence the free $A$-module of rank one detects weak equivalences, as claimed.
\end{example}

Now we formulate the precise criterion for a Quillen functor
pair to be weakly monoidal. The proof of this proposition is
in Section \ref{sec-proofs}.

\begin{proposition}\label{prop-generating criterion}
Consider a Quillen functor pair  
$(\!\xymatrix{\lambda\co \cD \ar@<-.4ex>[r]&\cC:R\ar@<-.4ex>[l]}\!)$ 
between monoidal model categories together with a lax monoidal structure
on the right adjoint $R$. Suppose further that
\begin{enumerate} 
\item  for some (hence any) cofibrant replacement
$q\co\mID^c\to\mID$ of the unit object in $\cD$, the composite map
\[ \lambda(\mID^c)\ \xrightarrow{\, \lambda(q)\, } \
\lambda(\mID) \ \xrightarrow{\ \tilde\nu\ } \ \mIC \] 
is a weak equivalence in $\cC$,
where the second map is adjoint to the monoidal structure map 
$\nu\co \mID\to R(\mIC)$;
\item the unit object $\mID$ detects weak equivalences in $\cD$.
\end{enumerate}
Then $R$ and $\lambda$ are a weak monoidal Quillen pair.
\end{proposition}

This criterion works well in unstable situations such as (non-negatively
graded) chain complexes and simplicial abelian groups.  For the  
stable case, though, this notion of detecting weak equivalences is
often too strong.  Thus, we say an object $A$ of a stable model category
$\cC$ {\em stably detects weak equivalences} if $f\co Y \to Z$ is
a weak equivalence if and only if $[A, Y]_*^{\Ho\cC} \to [A, Z]_*^{\Ho\cC}$
is an isomorphism of the $\mZ$-graded abelian groups of morphisms in  
the triangulated homotopy category $\Ho \cC$.  For example, a (weak) generator
stably detects weak equivalences by~\cite[2.2.1]{stable}.   The unstable
notion above would correspond to just considering $*\geq 0$.     

Modifying the proof of Proposition~\ref{prop-generating criterion} 
by using the graded morphisms in the homotopy category instead
of the mapping spaces introduced above proves the following stable
criterion.

\begin{proposition}\label{prop-stable generating criterion}
Consider a Quillen functor pair  
$(\!\xymatrix{\lambda\co \cD \ar@<-.4ex>[r]&\cC:R\ar@<-.4ex>[l]}\!)$ 
between monoidal stable model categories together with a lax monoidal structure
on the right adjoint $R$. Suppose further that
\begin{enumerate} 
\item  for some (hence any) cofibrant replacement
$q\co\mID^c\to\mID$ of the unit object in $\cD$, the composite map
\[ \lambda(\mID^c)\ \xrightarrow{\, \lambda(q)\, } \
\lambda(\mID) \ \xrightarrow{\ \tilde\nu\ } \ \mIC \] 
is a weak equivalence in $\cC$,
where the second map is adjoint to the monoidal structure map 
$\nu\co \mID\to R(\mIC)$;
\item the unit object $\mID$ stably detects weak equivalences in $\cD$.
\end{enumerate}
Then $R$ and $\lambda$ are a weak monoidal Quillen pair.
\end{proposition}

\section{Chain complexes and simplicial abelian groups, revisited}
\label{sec-sAb and ch revisited}

Dold and Kan showed that the category of non-negatively graded chain 
complexes is equivalent to the category of simplicial abelian 
groups, see for example \cite[Thm.~1.9]{dold}. 
The equivalence is given by the normalization functor $N\co \sab\to \ch$ 
and its inverse $\Gam\co \ch\to \sab$.  
Because the two functors are inverses to each other, 
they are also adjoint to each other on both sides.
Hence the normalized chain complex  $N$ and 
its inverse functor $\Gamma$ give rise to two different 
weak monoidal Quillen equivalences. 
Each choice of right adjoint comes with a particular 
monoidal transformation, namely the shuffle map
(for $N$ as the right adjoint, see Section~\ref{N as right adjoint})
or the Alexander-Whitney map
(for $\Gam$ as the right adjoint, see Section~\ref{Gamma as right adjoint}).

In Section~\ref{sub-model structures} we recall the
supporting model category structures on chain complexes, simplicial
abelian groups and the associated categories of monoids and modules.
In Section~\ref{N as right adjoint}, we show that
Theorem~\ref{thm-WMQE gives equivalent algebras} parts (\ref{TWO}) and 
(\ref{THREE}) imply Theorem~\ref{thm-sab-ch} parts (\ref{ii}) and (\ref{iii}) 
respectively.  In Section~\ref{Gamma as right adjoint}, we show that
Theorem~\ref{thm-WMQE gives equivalent algebras} part (\ref{TWO}) 
implies Theorem~\ref{thm-sab-ch} part (\ref{i}).  
In Section~\ref{algebras over com sRing} we then use the fact 
that the shuffle map for $N$ is lax symmetric monoidal to deduce 
Theorem~\ref{thm-sab-ch} part (\ref{iv}) from Theorem~\ref{thm-WMQE gives 
equivalent algebras} part (\ref{THREE}).

\subsection{Model structures} \label{sub-model structures}
Let $k$ be a commutative ring. 
The category $\ch_k$ of non-negatively graded chain complexes of $k$-modules
supports the {\em projective} model structure:
weak equivalences are the quasi-isomorphisms, fibrations are the chain 
maps which are surjective in positive dimensions, and cofibrations
are the monomorphisms with dimensionwise projective cokernel.
More details can be found in~\cite[Sec.~7]{ds} 
or \cite[2.3.11, 4.2.13]{hovey-book}
(the references in~\cite{hovey-book} actually treat $\mZ$-graded chain
complexes, but the arguments for $\ch$ are similar; there is also an
{\em injective} model structure for $\mZ$-graded chain complexes 
with the same weak equivalences, but we do not use this model structure).

The model category structure on simplicial $k$-modules has as weak equivalences
and fibrations the weak homotopy equivalences and Kan fibrations
on underlying simplicial sets; the cofibrations are the retracts of 
free maps in the sense of~\cite[II.4.11 Rem.\ 4]{Q}.
For more details see ~\cite[II.4, II.6]{Q}.

The functor $N$ is an inverse equivalence of categories,
and it exactly matches the notions of cofibrations, fibrations and
weak equivalences in the two above model structures. 
So $N$ and its inverse $\Gamma$ can be viewed as Quillen equivalences 
in two ways, with either functor playing the role of the left 
or the right adjoint.

Both model structures are compatible with the tensor products
of Section \ref{sub-Tensor products} in the sense that they satisfy
the pushout product axiom (compare~\cite[4.2.13]{hovey-book}).
A model category on simplicial rings with fibrations and weak equivalences
determined on the underlying simplicial abelian groups (which in
turn are determined by the underlying simplicial sets) was established
by Quillen in~\cite[II.4, Theorem 4]{Q}.  Similarly, there is a model structure
on differential graded rings with weak equivalences the quasi-isomorphisms
and fibrations the maps which are surjective in positive degrees;
see~\cite{jardine}. These model category structures also follow from
verifying the monoid axiom~\cite[Def.~3.3]{ss} and using 
Theorem 4.1 of~\cite{ss}; see also~\cite[Section 5]{ss}.

\subsection{A first weak monoidal Quillen equivalence}
\label{N as right adjoint} 
Let $k$ be a commutative ring. 
We view the normalization functor 
\[ N \  : \  s\Modr k\ \to\ \ch_k \] 
from simplicial $k$-modules
to non-negatively graded chain complexes of $k$-modules as the 
{\em right adjoint} and its inverse $\Gamma$ \eqref{eq-def Gamma}
as the left adjoint.

We consider $N$ as a lax monoidal functor via the shuffle map
\eqref{def-restricted shuffle}.
The shuffle map $\nabla\co NA \tensor NB\to N(A\tensor B)$
is a chain homotopy equivalence for every pair of simplicial $k$-modules, 
cofibrant or not, with homotopy inverse 
the Alexander-Whitney map~\cite[29.10]{may}. 
Since $\Gamma$ takes quasi-isomorphisms 
to weak equivalences of simplicial $k$-modules,
and since the unit and counit of the adjunction between $N$ and $\Gamma$
are isomorphisms, the description \eqref{eq-tvarphi} of
the comonoidal transformation for the left adjoint shows that the map
\[ \widetilde\nabla_{C,D} \ : \ \Gamma(C \tensor D) \ \to \ 
\Gamma C\tensor \Gamma D  \] 
is a weak equivalence of simplicial $k$-modules for all 
connective complexes of $k$-modules.

In other words, with respect to the shuffle map,
$N$ is the right adjoint of a weak monoidal Quillen equivalence
between simplicial $k$-modules and connective chain complexes of $k$-modules.
Since the unit objects are cofibrant, we can apply 
Theorem \ref{thm-WMQE gives equivalent algebras}.
Part \eqref{THREE} shows that normalization is the right adjoint
of a Quillen equivalence from simplicial $k$-algebras to connective
differential graded $k$-algebras; 
this proves part \eqref{iii} of Theorem~\ref{thm-sab-ch}.
Quillen invariance holds for simplicial rings, and normalization
preserves all weak equivalences; so for $k=\mZ$, part \eqref{TWO} 
of Theorem \ref{thm-WMQE gives equivalent algebras} 
shows that for every simplicial ring $A$, normalization is the right adjoint
of a Quillen equivalence from simplicial $A$-modules to connective
differential graded $NA$-modules;  this proves part \eqref{ii}
of Theorem~\ref{thm-sab-ch}.

All the above does not use the fact that the shuffle map 
for $N$ is {\em symmetric} monoidal.
We explore the consequences of this in \ref{algebras over com sRing}.

\subsection{Another weak monoidal Quillen equivalence}
\label{Gamma as right adjoint}
Again, let $k$ be a commutative ring. 
This time we treat the normalization functor 
$N\co s\Modr k\to \ch_k$ as the {\em left} adjoint and 
its inverse functor $\Gamma\co \ch_k\to s\Modr k$
as {\em right} adjoint to $N$.
The monoidal structure on $\Gam$ defined in \eqref{def-AW as monoidal} 
is made so that the comonoidal transformation \eqref{eq-tvarphi}
for the left adjoint $N$ is precisely the Alexander-Whitney map
\[ AW \ : \ N(A \tensor B) \ \to \ NA\tensor NB \ .  \] 
The Alexander-Whitney map is a chain homotopy equivalence for arbitrary
simplicial $k$-modules $A$ and $B$, with homotopy inverse 
the shuffle map~\cite[29.10]{may}; hence $\Gamma$ becomes the
right adjoint of a weak monoidal Quillen equivalence.

Since the unit objects are cofibrant, we can again apply 
Theorem \ref{thm-WMQE gives equivalent algebras}.
Part \eqref{TWO} shows that for every connective differential graded ring $R$,
the functor $\Gamma$ is the right adjoint of a Quillen equivalence
from connective differential graded $R$-modules 
to simplicial $\Gamma R$-modules;  
this proves part \eqref{i} of Theorem~\ref{thm-sab-ch}.

\subsection{Modules and algebras over a commutative simplicial ring}
\label{algebras over com sRing}
The shuffle map \eqref{def-shuffle} is lax {\em symmetric} monoidal,
and so is its extension
\[ \nabla \ : \ NA \tensor NB \ \to \ N(A\tensor B) \]
to normalized chain complexes.
In sharp contrast to this, the Alexander-Whitney map is {\em not} symmetric.
This has the following consequences:

$\bullet$\qua If $A$ is a {\em commutative} simplicial ring, 
then the normalized chains
$NA$ form a differential graded algebra which is commutative 
in the graded sense, i.e., we have
\[ xy \ = \ (-1)^{|x|\,|y|} yx \] 
for homogeneous elements $x$ and $y$ in $NA$.

$\bullet$\qua The functor $N$ inherits a lax monoidal structure 
when considered as a functor from simplicial $A$-modules 
(with tensor product over $A$)
to connective differential graded $NA$-modules 
(with tensor product over $NA$). More precisely, there is a 
unique natural chain map
\[ \nabla^A \ : \ NM \tensor_{NA} N(M') \ \to \  N(M\tensor_A M') \ . \] 
for $A$-modules $M$ and $M'$, such that the square 
\[\xymatrix@C=20mm{NM \tensor NM' \ar[r]^{\nabla} \ar[d] & 
N(M \tensor M') \ar[d] \\
NM\tensor_{NA}\tensor NM' \ar[r]_{\nabla^A} & N(M\tensor_A M') }\]
commutes (where the vertical morphisms are the natural quotient maps).
This much does not depend on commutativity of the product of $A$. 
However, if $A$ is commutative, then $\nabla^A$ constitutes a 
lax symmetric monoidal functor from $A$-modules to $NA$-modules;
this uses implicitly that the monoidal transformation $\nabla$
is symmetric, hence compatible with the isomorphism 
of categories between left and right modules over $A$ and $NA$. 

Now let $L^A\co \Modr NA\to \Modr A$ denote the left adjoint of $N$ 
when viewed as a functor from left $A$-modules to left modules over $NA$
(compare Section \ref{sub-left adjoint}).
Then the lax comonoidal map for $L^A$ has the form
\[ L^A(W\tensor_{NA} W') \ \to \ L^A(W)\tensor_{A}L^A(W') \] 
for a pair of left $NA$-modules $W$ and $W'$.
We claim that this map is a weak equivalence for cofibrant modules
$W$ and $W'$ by appealing to Theorem \ref{prop-generating criterion}.
Indeed, the unit objects of the two tensor products are the respective free
modules of rank one, which are cofibrant.
The unit map $\eta \co L^A(NA) \to A$ is even an isomorphism,
and the free $NA$-module of rank one detects weak equivalences 
by Example \ref{ex-free rank one generates}.
So Theorem \ref{prop-generating criterion} applies
to show that the adjoint functor pair
\[\xymatrix@=12mm{ \Modr NA \quad \ar@<-.4ex>_-{L^A}[r] &
\quad  \Modr A\quad  \ar@<-.4ex>_-{N}[l] }   \]
is a weak monoidal Quillen pair.
These two adjoint functors form a Quillen equivalence by 
part \eqref{TWO} of Theorem~\ref{thm-WMQE gives equivalent algebras}.

Now we can apply part \eqref{THREE} of
Theorem \ref{thm-WMQE gives equivalent algebras}.
The conclusion is that the normalized chain complex functor
is the right adjoint of a Quillen equivalence from the
model category of simplicial $A$-algebras to the model category
of connective differential graded algebras over the
commutative differential graded ring $NA$. 
In other words, this proves part \eqref{iv} of Theorem \ref{thm-sab-ch}

\section{Proofs}\label{sec-proofs}

This section contains the proofs of the main results of the paper,
namely Theorem \ref{thm-WMQE gives equivalent algebras}
and the criterion for being weakly monoidal, 
Proposition \ref{prop-generating criterion}.

\subsection{Proof of Theorem \ref{thm-WMQE gives equivalent algebras}}
This proof depends on a comparison of different kinds of left adjoints
in Proposition \ref{prop-L}.
For this part we assume that $\lambda$ and $R$ 
form a weak monoidal Quillen pair, 
in the sense of Definition \ref{def-WMQ equivalence}, 
between two monoidal model categories 
$(\cC, \otimes, \mIC)$ and  $(\cD, \sm, \mID )$.
As before, the lax monoidal transformation of $R$ is denoted 
$\varphi_{X,Y}\co  RX \sm RY\to R(X\otimes Y)$;
it is `adjoint' to  the comonoidal transformation \eqref{eq-tvarphi}
\[ \tvarphi \ \co \ \lambda( A \sm B)\ \to \ \lambda A \otimes \lambda B 
\ . \]
This comonoidal map need not be an isomorphism; in that case $\lambda$ does not
have a {\em monoidal} structure, and so it does not pass to a functor
between monoid and module categories.
However, part of the definition of a weak monoidal Quillen pair
is that $\tvarphi$ is a weak equivalence whenever $A$ and $B$ are cofibrant.

In our applications, the monoidal functor $R$ has left adjoints
on the level of monoids and  modules (see Section \ref{sub-left adjoint}),
and the following proposition compares these `structured' left adjoints 
to the underlying left adjoint $\lambda$.
If the monoidal transformation $\tvarphi$ of~\eqref{eq-tvarphi}
happens to be an isomorphism, then so are the comparison morphisms
$\chi_B$ and $\chi_W$ which occur in the following proposition; 
so for {\em strong} monoidal Quillen pairs
the following proposition has no content.

In Definition \ref{def-create} we defined the notion of a `cell object'
relative to a functor to a model category.
In the following proposition, the notions of cell object are taken
relative to the forgetful functors from algebras, respectively
modules, to the underlying monoidal model category.

\begin{proposition}\label{prop-L}
Let  $(\!\xymatrix{\lambda\co \cD \ar@<-.4ex>[r]&\cC:R\ar@<-.4ex>[l]}\!)$ 
be a weak monoidal Quillen pair between monoidal model categories 
with cofibrant unit objects.
\begin{enumerate}
\item  Suppose that the functor $R\co \cC\mon\to\cD\mon$
has a left adjoint $L^{\text{\em mon}}$.
Then for every cell monoid $B$ in $\cD$, the $\cC$-morphism
\[ \chi_B \ \co \ \lambda B \ \to \ L^{\text{\em mon}}B \]  
which is adjoint to the underlying $\cD$-morphism of the adjunction unit
$B\to R(L^{\text{\em mon}}B)$ is a weak equivalence.
\item Let $B$ be a cell monoid in $\cD$
for which the functor $R\co \Modr(L^{\text{\em mon}}B)\to\Modr B$
has a left adjoint $L_B$.
Then for every cell $B$-module $M$, the $\cC$-morphism 
\[ \chi_M \ \co \ \lambda M \ \to \ L_BM \]  
which is adjoint to the underlying $\cD$-morphism of the adjunction unit
$M\to R(L_BM)$ is a weak equivalence.
\end{enumerate}
\end{proposition}
\begin{proof} Part (1): The left adjoint $L^{\text{mon}}$ takes the initial
$\cD$-monoid $\mID$ to the initial $\cC$-monoid $\mIC$ and the map
$\chi_{\mID}\co \lambda(\mID)\to L^{\text{mon}}\mID\iso \mIC$ 
is the adjoint $\tilde\nu$ of the unit map $\nu\co \mID\to R(\mIC)$.
By definition of a weak monoidal Quillen pair, $\tilde\nu$
is a weak equivalence. Hence $\chi$ is a weak equivalence for the
initial $\cD$-monoid $\mID$.   

Now we proceed by a cell induction argument, along free extensions of
cofibrations in $\cD$.
We assume that $\chi_B$ is a weak equivalence 
for some cell $\cD$-monoid $B$. Since the unit $\mID$ is cofibrant, 
$B$ is also cofibrant in the underlying category $\cD$ 
by an inductive application of \cite[Lemma 6.2]{ss}.

We consider another monoid $P$ obtained from $B$ by a single 
`cell attachment', i.e., a pushout in the category 
of $\cD$-monoids of the form
\[\xymatrix{ T_{\cD} K \ar[r] \ar[d] &
T_{\cD} K' \ar[d] \\ B \ar[r] & P }\]
where $K \to K'$ is a cofibration in $\cD$. We may assume  without loss
of generality that $K$ and $K'$ are in fact cofibrant. 
Indeed, $P$ is also a pushout of the diagram
\[\xymatrix{ B & T_{\cD} B \ar[r] \ar[l]_{\epsilon} &T_{\cD}(B\amalg_K K') }\]
where $\epsilon\co  T_{\cD}B\to B$ is the counit of the free monoid adjunction
and $B\amalg_KK'$ denotes the pushout in the underlying category $\cD$.
Since $B$ is cofibrant in $\cD$, the morphism $B\to B\amalg_K K'$ 
is a cofibration between cofibrant objects in $\cD$, 
and it can be used instead of the original cofibration $K\to K'$.

Free extensions of monoids are analyzed 
in the proof of~\cite[6.2]{ss} and we make use of that description. 
The underlying object of $P$ can be written as a colimit 
of a sequence of cofibrations $P_{n-1} \to P_n$  in $\cD$,
with $P_0=B$ such that each morphism $P_{n-1} \to P_n$ 
is a pushout in $\cD$ of a particular cofibration 
$Q_n(K,K',B) \to (B \sm K')^{\sm n} \sm B$.
Since $\lambda$ is a left adjoint on these underlying categories, $\lambda$ 
applied to each corner of these pushouts is still a pushout square. 

Since $L^{\text{mon}}$ is a left adjoint on the categories of monoids, 
it preserves pushouts of monoids.  
So we have the following pushout of $\cC$-monoids
\[\xymatrix{ L^{\text{mon}}(T_{\cD}K) \ar[r] \ar[d] &
L^{\text{mon}}(T_{\cD}K') \ar[d] \\
L^{\text{mon}}B \ar[r] & L^{\text{mon}}P }\]
Because of the natural isomorphism \eqref{eq-free monoid iso} between
$L^{\text{mon}}(T_{\cD} K)$ and $T_{\cC}(\lambda K)$,
the pushout $L^{\text{mon}}P$ is thus the free extension 
of $L^{\text{mon}}B$ along the cofibration $\lambda K\to\lambda K'$
between cofibrant objects in $\cC$.
In particular, $L^{\text{mon}}P$ is the colimit in $\cC$ 
of a sequence of cofibrations $R_{n-1} \to R_n$ 
each of which is a pushout in $\cC$ of a cofibration
$Q_n(\lambda K, \lambda K', L^{\text{mon}}B) \to 
(L^{\text{mon}}B \sm \lambda K')^{\sm n} \sm L^{\text{mon}}B$.

The map $\chi_P\co \lambda P \to L^{\text{mon}}P$ 
preserves these filtrations: it takes $\lambda P_n$ to $R_n$.  
We now show by induction that for each $n$ the
map $\lambda P_n \to R_n$ is a weak equivalence in $\cC$.  We show that
the map on each of the other three corners of the pushout squares defining 
$\lambda P_n$ and $R_n$ is a weak equivalence between cofibrant objects;
then we apply~\cite[5.2.6]{hovey-book}
to conclude that the map of pushouts is also a weak equivalence.
By induction we assume $\lambda P_{n-1} \to R_{n-1}$ is a weak equivalence.
A second corner factors as 
\begin{align*} \lambda ((B \sm K')^{\sm n} \sm B) \ 
&\xrightarrow{\qquad\ \, \tvarphi\qquad\ } \qquad
(\lambda B \sm \lambda K')^{\sm n} \sm \lambda B \\ 
&\xrightarrow{(\chi_B\sm\Id)^{\sm n}\sm\chi_B} \ 
(L^{\text{mon}}B \sm \lambda K')^{\sm n} \sm L^{\text{mon}}B \ . 
\end{align*}
Since the first map is an (iterated) instance
of the comonoidal transformation $\tvarphi$, and since $B$ is cofibrant
in the underlying category $\cD$, the first map is a weak equivalence
by hypothesis. Since $L^{\text{mon}}B$ is a cell monoid in $\cC$ 
and the unit object $\mIC$ is cofibrant, $L^{\text{mon}}B$ is
cofibrant in the underlying category $\cC$ by~\cite[Lemma 6.2]{ss}.
By induction we know that $\chi_B\co \lambda B\to L^{\text{mon}}B$ 
is a weak equivalence;
since all objects in sight are cofibrant and smashing with a cofibrant object
preserves weak equivalences between cofibrant objects, the second map is 
also a weak equivalence.  

The third corner works similarly: the map factors as 
\[ \lambda Q_n(K, K', B) \ \to \ Q_n(\lambda K, \lambda K', \lambda B) 
\ \to \ Q_n(\lambda K, \lambda K', L^{\text{mon}}B) \ . \]  
Here $Q_n$ itself is constructed as a pushout of an $n$-cube 
where each map is a cofibration between cofibrant objects.  
Using a variant of~\cite[5.2.6]{hovey-book} for $n$-cubes, 
the hypothesis on the comonoidal transformation $\tvarphi$ 
and induction shows that this third corner is also a weak equivalence.   
	   
Since each filtration map is a cofibration between cofibrant objects
and $\lambda P_n \to R_n$ is a weak equivalence for each $n$, we can 
apply~\cite[1.1.12, 5.1.5]{hovey-book} to conclude that the map of 
colimits $\lambda P \to L^{\text{mon}}P$ is a weak equivalence since it
is the colimit of a weak equivalence between cofibrant objects in the
Reedy model category of directed diagrams.  
Similarly, for the transfinite compositions allowed in building
up a cell object,~\cite[1.1.12, 5.1.5]{hovey-book} shows that the
map of colimits is a weak equivalence. 

\medskip

(2)\qua We use a similar induction for cell $B$-modules.
Again at transfinite composition steps, 
\cite[1.1.12, 5.1.5]{hovey-book} gives the necessary conclusion;
so we are left with considering the single cell attachments.
Suppose the statement has been verified for some cell $B$-module $M$. 
Since $B$ is a cell monoid, it is cofibrant in the underlying category $\cD$,
and so is $M$, by induction on the `number of cells'.
Suppose $M'$ is obtained from $M$ as a free extension of $B$-modules, i.e,
it is a pushout of a diagram
\[\xymatrix {M & K\sm B \ar[l] \ar[r] & K'\sm B }\]
Where $K \to K'$ is a cofibration in $\cD$.
By the same trick as in the first part, we can assume 
without loss of generality that $K$ and $K'$ are cofibrant in $\cD$; 
this exploits the fact that $B$ is also cofibrant in $\cD$.  

Since $L_B$ is a left adjoint, it preserves pushouts.
On the other hand, the forgetful functor from $B$-modules to $\cD$
and $\lambda\co \cC\to\cD$ also preserve pushouts, so we get a 
commutative diagram in $\cC$
\[\xymatrix@C=4mm{ \lambda(K\sm B) \ar[dd]\ar[rrr]^-{\chi_{K\sm B}}\ar[dr] &&& 
L_B(K\sm B) \ar'[d][dd] \ar[dr] \\
& \lambda(K'\sm B) \ar[rrr]^(.3){\chi_{K'\sm B}} \ar[dd] &&& 
L_B(K'\sm B) \ar[dd] \\
\lambda M \ar'[r][rrr]^{\chi_M} \ar[dr] &&& L_BM \ar[dr] \\
& \lambda M' \ar[rrr]_{\chi_{M'}}  &&& L_BM'
}\]
in which the left and right faces are pushout squares.
Because of the natural isomorphisms 
$L_B(K\sm B)\iso (\lambda A)\tensor L^{\text{mon}}B$ 
and $L_B(K'\sm B)\iso \lambda K'\tensor L^{\text{mon}}B$
of \eqref{eq-free module iso for L_B}, the pushout $L_BM'$ is thus 
the free extension of $L_BM$ along the cofibration $\lambda K\to\lambda K'$
between cofibrant objects in $\cC$.

By assumption, the map $\chi_M$ is a weak equivalence.
The map $\chi_{K\sm B}$  factors as the composite 
\begin{equation}\label{composite} 
\lambda (K \sm B) \ \xrightarrow{\ \tvarphi\ } \
\lambda K \tensor \lambda B \ \xrightarrow{\text{Id}\tensor\chi_B} \ 
\lambda K \tensor L^{\text{mon}}B\ . \end{equation}
Since $B$ is cofibrant in the underlying category $\cD$
and the first map is an instance of the comonoidal transformation $\tvarphi$,
it is a weak equivalence.
Since the map $\lambda B\to L^{\text{mon}}B$ is a weak equivalence 
by Part (1) and all objects in sight are cofibrant, 
the second map of \eqref{composite} is also a weak equivalence.
Hence $\chi_{K\sm B}$ is a weak equivalence, 
and similarly for $\chi_{K'\sm B}$.
 Since the maps on each of the three initial corners 
of the pushout squares defining $\lambda M'$ and $L_BM'$ are 
weak equivalence between cofibrant objects,
the map of pushouts is also a weak equivalence
(see for example~\cite[5.2.6]{hovey-book}).
So $\chi_{M'}\co \lambda M' \to L_BM'$ is a weak equivalence.
\end{proof}

\begin{proof}[Proof of Theorem~\ref{thm-WMQE gives equivalent algebras}]
Since the fibrations and trivial fibrations of monoids and modules
are defined on the underlying category, 
the right adjoint $R$ is a right Quillen functor in all cases.

\eqref{THREE}\qua Consider a cofibrant $\cD$-monoid $B$ 
and a fibrant $\cC$-monoid $Y$. We have to show that a monoid homomorphism
$B\to RY$ is a weak equivalence 
if and only if its adjoint $L^{\text{mon}}B \to Y$ is a weak equivalence.
In the underlying category $\cC$, we can consider the composite
\begin{equation}\label{monoid vs underlying adjoint} 
\lambda B \ \xrightarrow{\chi_B} \ L^{\text{mon}}B \ \to \ Y \end{equation}
which is adjoint to the underlying $\cD$-morphism of $B\to RY$.

Since the forgetful functor creates 
(in the sense of Definition \ref{def-create}) the model
structure in the category of $\cD$-monoids, every cofibrant $\cD$-monoid $B$
is a retract of a cell monoid by definition.
So the morphism $\chi_B\co \lambda B \to L^{\text{mon}}B$ is a weak equivalence 
by Proposition~\ref{prop-L} (1).
Since $\mID$ is cofibrant, a cofibrant $\cD$-monoid is also cofibrant as an 
object of $\cD$ by~\cite[6.2]{ss}.  Thus, since $\lambda$ and $R$ form 
a Quillen equivalence on the underlying categories, 
the composite  map \eqref{monoid vs underlying adjoint}
is a weak equivalence if and only if $B \to RY$ is a weak equivalence.  
So $L^{\text{mon}}$ and $R$ form a Quillen equivalence between the
categories of monoids in $\cC$ and $\cD$.

\eqref{ONE}\qua This is very similar to part~\eqref{THREE}, but using the
second part of Proposition~\ref{prop-L} instead of the first part.

\eqref{TWO}\qua Let $\psi\co A \to A^{\text{f}}$ be a fibrant
replacement in the category of $\cC$-monoids.
If $A$ is already fibrant, take $A^{\text{f}}=A$,
otherwise assume $R$ preserves all weak equivalences, so either way
$R\psi\co RA\to R(A^{\text{f}})$ is a weak equivalence
of $\cD$-monoids. Let $\kappa\co R(A^{\text{f}})^c\to R(A^{\text{f}})$ 
be a cofibrant replacement in the category of $\cD$-monoids.  
Let $\widetilde{\kappa} \co L^{\text{mon}}(R(A^{\text{f}})^c)\to A^{\text{f}}$
be its adjoint; by part \eqref{THREE}, this adjoint $\widetilde{\kappa}$ 
is a weak equivalence of monoids in $\cC$.
We have a commutative diagram of right Quillen functors
\[\xymatrix@C=12mm{ \Modr A \ar[d]_-{R} & 
\Modr A^{\text{f}} \ar[d]^-{R} \ar[l]_-{\psi^*}
 \ar[r]^-{\widetilde{\kappa}^*} &
\Modr L^{\text{mon}}(R(A^{\text{f}})^c) \ar[d]^{R} \\
\Modr RA & 
\Modr R(A^{\text{f}}) \ar[l]^-{(R\psi)^*} \ar[r]_{\kappa^*} &  
\Modr R(A^{\text{f}})^c}\]
in which the horizontal functors are restrictions of scalars along 
the various weak equivalences of monoids; 
these are right Quillen equivalences by Quillen invariance.
By part~\eqref{ONE} the right vertical functor 
is a right Quillen equivalence. Hence, the middle and left vertical functors 
are also right Quillen equivalences.
\end{proof}

\subsection{Proof of Proposition \ref{prop-generating criterion}}
We are given a Quillen functor pair  
$(\!\xymatrix{R\co \cC \ar@<-.4ex>[r]&\cD:\lambda\ar@<-.4ex>[l]}\!)$ 
between monoidal model categories together with a lax monoidal structure
on the right adjoint $R$. Moreover, the unit object $\mID$ detects 
weak equivalences in $\cD$ (in the sense of Definition \ref{def-detect})
and for some (hence any) cofibrant replacement
$\mID^c\to\mID$ of the unit object in $\cD$, the composite map
$\lambda(\mID^c)\to\lambda(\mID)\xrightarrow{\tilde\nu}\mIC $ 
is a weak equivalence in $\cC$.
We have to show that $R$ and $\lambda$ are a weak monoidal Quillen pair.

Our first step is to show that 
for every cofibrant object $B$ of $\cD$ and every fibrant object $Y$ of $\cC$,
a certain map 
\begin{equation}\label{eq-adjoint form}
 R\Hom_{\cC}(\lambda B,Y) \ \to \ \Hom_{\cD}(B,RY) \end{equation}
is a weak equivalence in $\cD$. Here $\Hom_{\cC}(-,-)$ and $\Hom_{\cD}(-,-)$
denote the internal function objects in $\cC$ respectively $\cD$,
which are part of the {\em closed} symmetric monoidal structures. 
The map \eqref{eq-adjoint form} is adjoint to the composition
\begin{align*} 
R\Hom_{\cC}(\lambda B,Y) \sm B \ &\xrightarrow{\Id\sm\eta} \ 
R\Hom_{\cC}(\lambda B,Y) \sm R(\lambda B) \\ 
&\xrightarrow{\ \ \varphi\ \, } \ 
R(\Hom_{\cC}(\lambda B,Y) \tensor \lambda B) \
\xrightarrow{R(\text{eval})} \ RY \ . \end{align*} 

Choose a cosimplicial resolution $\mID^*$ of the unit object; in particular,
in cosimplicial dimension zero we get a cofibrant replacement 
$\mID^0\stackrel{\sim}{\to}\mID$ of the unit object.
Since $\mID$ detects weak equivalences in $\cD$ 
and the objects $R\Hom_{\cC}(\lambda B,Y)$ and $\Hom_{\cD}(B,RY)$ are fibrant, 
we can prove that \eqref{eq-adjoint form} is a weak equivalence 
by showing that we get a weak equivalence of mapping spaces 
\[ \map_{\cD}(\mID^*, R\Hom_{\cC}(\lambda B,Y))\ \to\ 
\map_{\cD}(\mID^*,\Hom_{\cD}(B,RY))\ . \] 
The adjunction isomorphisms between $R$ and $\lambda$ and between
the monoidal products and internal function objects allow us to rewrite
this map as
\[ \map_{\cC}(\tvarphi,Y) \ : \
\map_{\cC}(\lambda(\mID^*)\tensor\lambda B,Y) \ \to \ 
\map_{\cC}(\lambda(\mID^*\sm B),Y) \ , \]
where the cosimplicial map
\begin{equation}\label{eq-cosimplicial varphi}
\tvarphi \ : \ \lambda(\mID^*\sm B)\ \to\ \lambda(\mID^*)\tensor\lambda B 
\end{equation}
is an instance of the comonoidal map $\tvarphi$ in each cosimplicial dimension.
We consider the commutative diagram in $\cC$
\[\xymatrix{\lambda(\mID^0\sm B)\ar[d]_{\sim} \ar[rr]^-{\tvarphi} &&
\lambda(\mID^0)\tensor\lambda B \ar[d]^{\sim} \\
\lambda(\mID\sm B) \ar[r]_-{\iso} & \ \lambda B \ & \mIC\tensor\lambda B 
\ar[l]^-{\iso} } \]
whose top horizontal map is the component 
of \eqref{eq-cosimplicial varphi} in cosimplicial dimension zero.
Since $\mID^0$ is a cofibrant replacement of the unit in $\cD$,
the map $\mID^0\sm B\to\mID \sm B\iso B$ is a weak equivalence 
between cofibrant
objects, and so the left vertical map is a weak equivalence.
By hypothesis (1), the composite map $\lambda(\mID^0)\to\lambda(\mID)\to\mIC $
is a cofibrant replacement of the unit in $\cC$, so
smashing it with the cofibrant object $\lambda B$ gives the right vertical
weak equivalence.

We conclude that \eqref{eq-cosimplicial varphi} is a weak equivalence 
in dimension zero; every left Quillen functor such as $\lambda$, $-\sm B$ 
or $-\tensor\lambda B$ preserves cosimplicial resolutions, 
so \eqref{eq-cosimplicial varphi} is a level equivalence 
between Reedy cofibrant objects. 
So the induced map on mapping spaces 
$\map(-,Y)$ is a weak equivalence by~\cite[5.4.8]{hovey-book},
and so the map \eqref{eq-adjoint form} is a weak equivalence in $\cD$.

Now we play the game backwards. If $A$ is another cofibrant object of $\cD$, 
then $\mID^*\sm A$ is a cosimplicial resolution whose dimension zero object
$\mID^0\sm A$ is weakly equivalent to $A$.
Since the map \eqref{eq-adjoint form} is a weak equivalence between fibrant
objects
and the functor $\map_{\cD}(\mID^*\sm A,-)$ 
is a right Quillen functor~\cite[5.4.8]{hovey-book}, 
we get an induced weak equivalence on mapping spaces
\[ \map_{\cD}(\mID^*\sm A, R\Hom_{\cC}(\lambda B,Y))\ \to\ 
\map_{\cD}(\mID^*\sm A,\Hom_{\cD}(B,RY))\ . \] 
Using adjunction isomorphisms again, this map can be rewritten as
\[ \map_{\cC}(\tvarphi,Y) \ : \ 
\map_{\cC}(\lambda(\mID^*\sm A) \tensor \lambda B,Y)\ \to\ 
\map_{\cC}(\lambda(\mID^*\sm A\sm B),Y)\ . \] 
where 
\[ \tvarphi \ : \ \lambda(\mID^*\sm A\sm B) \ \to \
\lambda(\mID^*\sm A)\tensor\lambda B \]
is another instance of the comonoidal map $\tvarphi$ 
in each cosimplicial dimension.

Since $\lambda(\mID^*\sm A)\tensor\lambda B$ is a cosimplicial frame of the
cofibrant object $\lambda A\tensor\lambda B$, 
the components of the mapping space
$\map_{\cC}(\lambda(\mID^*\sm A) \tensor \lambda B,Y)$ are in natural 
bijection with the morphisms from $\lambda A\tensor \lambda B$ 
in the homotopy category of $\cC$; similarly the components of
$\map_{\cC}(\lambda(\mID^*\sm A\sm B),Y)$ are isomorphic to
$[\lambda(A\sm B),Y]_{\Ho(\cC)}$. So we conclude that the comonoidal map 
$\tvarphi\co \lambda(A\sm B)\to \lambda A\tensor \lambda B$ induces a bijection
of homotopy classes of maps into  every fibrant object $Y$.
Thus the map $\tvarphi$ is a weak equivalence in $\cC$, 
and so we indeed have a weak monoidal Quillen pair.

\section{Some enriched model category theory}
\label{app}

In this section we develop some general theory of modules over 
`rings with several objects' based on a monoidal model category.
Most of this is a relatively straightforward generalization from
the case of `ring objects' or monoids to `rings with several objects'.
Throughout this section, $\cC$ is a monoidal model category 
(Definition~\ref{def-mmc}) with product $\tensor$ and unit object $\mI$.
For the special case where $\cC$ is the category of symmetric spectra,
this material can also be found in \cite[Sec.~A.1]{stable}.
Additional material on the homotopy theory of $\cC$-categories
can be found in~\cite{dundas}.

\subsection{Modules over $\cC$-categories}

Let $I$ be any set and let $\cO$ be a 
$\cC I$-category~\cite[6.2]{bor}, i.e., 
a category enriched over $\cC$ whose set of objects is $I$. 
This means that for all $i,j\in I$ 
there is a morphism object $\cO(i,j) \in \cC$,
unit morphisms $\mI\to \cO(i,i)$ and coherently associative 
and unital composition morphisms
\[ \cO(j,k) \tensor \cO(i,j) \ \to \ \cO(i,k) \ . \]
One may think of $\cO$ as a
`ring/monoid with many objects', indexed by $I$, enriched in $\cC$. 
Indeed if $I=\{*\}$ has only one element, then $\cO$ is completely determined
by the endomorphism $\cC$-monoid $\cO(*,*)$.
Moreover the $\cO$-modules as defined below coincide with the 
$\cO(*,*)$-modules in the ordinary sense.

A (right) {\em $\cO$-module} is a {\em contravariant} 
$\cC$-functor~\cite[6.2.3]{bor} from $\cO$ to the category $\cC$; 
explicitly, an  $\cO$-module $M$ consists of $\cC$-objects $M(i)$ 
for $i\in I$ and $\cC$-morphisms
\[ M(j) \tensor \cO(i,j) \ \to \ M(i) \] 
which are appropriately associative and unital. 
A map of $\cO$-modules is a $\cC$-natural transformation~\cite[6.2.4]{bor}.
A map of $\cO$-modules\  $f \co X \to Y$ is an
{\em objectwise equivalence} (or \ {\em objectwise fibration})
if $f(i)\co X(i) \to Y(i)$ is a weak equivalence
(fibration) in $\cC$ for each all $i\in I$.  A {\em cofibration}
is a map with the left lifting property with respect to any trivial fibration.
For every element $j\in I$, there is a {\em free} or {\em representable}
$\cO$-module $F_j^{\cO}$ defined by $F_j^{\cO}(i)=\cO(i,j)$.
As the name suggests, homomorphisms from $F_j^{\cO}$ into a module $M$
are in bijective correspondence with morphisms from $\mI$ to $M(j)$ 
by the enriched Yoneda Lemma~\cite[6.3.5]{bor}.  
The evaluation functor at $j$,
$\Ev_j \co \Modr\cO \to \cC$ has a left adjoint
`free' functor which sends an object $X$ of $\cC$ to the
module  $F_j^{\cO}X = X\tensor F_j^{\cO}$, where $X \tensor -$ is
defined by tensoring objectwise with $X$.

A {\em morphism} $\Psi\co \cO\to \cR$ of $\cC I$-categories is simply
a $\cC$-functor which is the identity on objects. 
We denote by $\cC I\cat$ the category of all $\cC I$-categories.
The {\em restriction of scalars}
\[ \Psi^* \ : \ \Modr\cR \ \to \ \Modr\cO \quad ; 
\quad M \ \longmapsto \ M \circ \Psi \]
has a left adjoint functor $\Psi_*$, also denoted $-\tensor_{\cO} \cR$,
which we refer to as {\em extension of scalars}. 
As usual it is given by an enriched coend,
i.e., for an $\cO$-module $M$ the $\cR$-module $\Psi_* M = M\tensor_{\cO} \cR$ 
is given by the coequalizer of the two $\cR$-module homomorphisms 
\[\xymatrix{ 
\bigvee_{i,j\in I} M(j) \, \tensor \, \cO(i,j)\, \tensor \, 
F^{\cR}_{i} \ \ar@<.3ex>[r] \ar@<-.3ex>[r] &          
\ \bigvee_{i\in I} M(i)\, \tensor\, F^{\cR}_{i} } \ . \] 
We call $\Psi\co \cO\to\cR$ a (pointwise) 
{\em weak equivalence} of $\cC I$-categories
if the $\cC$-morphism $\Psi_{i,j}\co \cO(i,j)\to\cR(i,j)$ 
is a weak equivalence for all $i,j\in I$, . 
Next we establish the model category structure
for $\cO$-modules and discuss Quillen invariance for
weak equivalences of $\cC$-categories.

\begin{theorem} \label{a.1.2}
Let $\cC$ be a cofibrantly generated monoidal model category 
which satisfies the monoid axiom~\cite[3.3]{ss} 
and such that every object of $\cC$ is small relative to the whole category. 
\begin{enumerate}
\item Let $\cO$ be a $\cC I$-category.
Then the category of $\cO$-modules with the objectwise equivalences,
objectwise fibrations, and cofibrations is a cofibrantly generated
model category. 
\item 
Let $\Psi\co \cO \to \cR$ be a weak equivalence of $\cC I$-categories.
Suppose that for every cofibrant right $\cO$-module $N$, 
the induced map $N \tensor_{\cO} \cO \to N \tensor_{\cO} \cR$ 
is an objectwise weak equivalence.  
Then restriction and extension of scalars along $\Psi$ form a $\cC$-Quillen
equivalence of the module categories.
\end{enumerate}
\end{theorem}
\begin{proof}
We use~\cite[2.3]{ss} to establish the model category for $\cO$-modules. 
Let $\mI_I$ denote the initial $\cC I$-category with $\mI_I(i,i)=\mI_{\cC}$ 
and $\mI_I(i,j)=\emptyset$, the initial object, for $i\ne j$. 
The category of $\mI_I$-modules is the $I$-indexed product category 
of copies of $\cC$. Hence it has a cofibrantly generated model category
inherited from $\cC$ in which the cofibrations, fibrations and weak
equivalences are objectwise.  Here the generating trivial cofibrations 
are the generating trivial cofibrations of $\cC$ between objects 
concentrated at one object, i.e. of the form $A_j$
with $A_j(j)= A$ and $A_j(i)=\emptyset$ if $i \ne j$.  

The unit morphism $\mI_I\to \cO$ induces adjoint functors of
restriction and extension of scalars between the module categories.  
This produces a triple 
$- \tensor_{\mI_I} \cO$ on $\mI_I$-modules with the algebras over this triple
the $\cO$-modules.  
Then the generating trivial cofibrations for $\cO$-modules are maps between 
modules of the form $A_j \tensor_{\mI_I} \cO = A \tensor \cO(-,j)$.    
Hence the monoid axiom applies to show 
that the new generating trivial cofibrations and their
relative cell morphisms are weak equivalences.  
Thus, since all objects in $\cC$ are assumed to be small, the model category
structure follows by criterion (1) of~\cite[2.3]{ss}.  

The proof of Part (2) follows as in~\cite[4.3]{ss}. 
The restriction functor $\Psi^*$ preserves objectwise fibrations and
objectwise equivalences, so restriction and extension of
scalars form a Quillen adjoint pair.
By assumption, for $N$ a cofibrant right $\cO$-module
 \[ N \ \iso \ N \tensor_{\cO} \cO \ \to \ N \tensor_{\cO} \cR \]
is a weak equivalence. Thus if $M$ is a fibrant right $\cR$-module, an
$\cO$-module map $N \to \Psi^* M$ is a weak equivalence if and only if the
adjoint $\cR$-module map $\Psi_* N = N \tensor_{\cO} \cR \to M$ is a weak 
equivalence.
\end{proof}

\subsection{Categories as monoids of graphs}
In~\cite[II.7]{ML-working} Mac Lane explains how a small category with object
set $I$ can be viewed as a monoid in the category of $I$-indexed graphs.
We need an enriched version of this giving $\cC$-categories as the monoids 
with respect to a monoidal product on the category of $\cC$-graphs, 
so that we can apply the general theory of \cite{ss}.  Note that here 
the product on $\cC$-graphs is {\em not} symmetric monoidal, so we must take
care in applying~\cite{ss}.

Let $(\cC,\tensor,\mIC)$ be a closed symmetric monoidal closed category
with an initial object $\emptyset$. Let $I$ be any set.
The category of (directed) {\em $I$-graphs in $\cC$}, 
denoted $\cC I\graph$ is simply the
product category of copies of $\cC$ indexed by the set $I\times I$.
If $G$ is an $I$-graph in $\cC$, then one can think of  
$G(i,j)$ as the $\cC$-object of arrows pointing from the vertex $i$ to the
vertex $j$. 

If $G$ and $H$ are two $I$-graphs, then their tensor product
is defined by the formula
\begin{equation}\label{eq-graph smash} 
(G\tensor H)(i,j) \ = \ \bigvee_{k\in I} \, G(k,j)\tensor H(i,k) \ . \end{equation}
The $I$-graph $\mI_I$ is defined by
\[ \mI_I(i,j) \ = \ \left\lbrace \begin{array}{ll} \mI & \mbox{ if } i=j \\
\emptyset & \mbox{ if } i\ne j \end{array} \right. \]
The smash product makes the category of $I$-graphs in $\cC$
into a monoidal category with unit object $\mI_I$ (but {\em not} a
symmetric monoidal category).
Moreover, the category $\cC I\cat$ of $\cC I$-categories is precisely 
the category of monoids in $\cC I\graph$ with respect to the smash product.  
Note that when $I$ is a singleton set $\cC I\graph$ is $\cC$ 
and $\cC I\cat$ is $\cC\mon$.

\medskip

{\bf Warning}\qua There is a slight risk of confusion in the notion of a
module over a $\cC I$-category when $I$ has more than one element.
As we just explained, such a $\cC I$-category $\cO$ 
is a monoid with respect to the monoidal product \eqref{eq-graph smash}
of $I$-graphs. So there is a notion of $\cO$-module which is an
$I$-graph $M$ together with a morphism of $I$-graphs $M\tensor \cO\to M$
which satisfies associativity and unit constraints. However, this is
{\em not} the same as the $\cO$-modules defined above as the enriched
functors from $\cO$ to $\cC$. These enriched functors have 
an underlying $I$-indexed family of objects in $\cC$, whereas
the other kind of modules have an underlying $I$-graph, 
so they have underlying $\cC$-objects indexed by {\em ordered pairs} of
elements from $I$. However, if $M$ is an $I$-graph with an associative,
unital $I$-graph morphism $M\tensor \cO\to M$, then we can fix an element 
$j\in I$, and obtain an enriched functor $M(-,j)$. So the modules which have
underlying $I$-graphs give rise to an $I$-indexed family of modules
in the earlier sense.

\medskip

Since $\cC I\graph$ is not a {\em symmetric} monoidal category, 
the results of~\cite{ss} do not apply directly
to produce a model category on the monoids, $\cC I\cat$.  
It turns out though that the proof Theorem 4.1 of~\cite{ss} carries over
since the homotopy type of a graph is determined pointwise 
and $\cC$ is assumed to be symmetric monoidal.   
First, if $\cC$ is a cofibrantly generated model category,
then $\cC I\graph$ is also a cofibrantly generated model category
with the cofibrations, fibrations and weak equivalences defined
{\em pointwise}, i.e., for each $(i,j)$.  The generating (trivial)
cofibrations are of the form $A_{i,j} \to B_{i,j}$ where $A \to B$
is a generating (trivial) cofibration in $\cC$ and $A_{i,j}$ is
the $\cC I\graph$ with value $A$ concentrated at $(i,j)$.  Based on this
model category we use \cite[4.1 (3)]{ss} to construct a model
category on the associated category of monoids, $\cC I\cat$.   

Dundas \cite[Thm.\ 3.3]{dundas} obtains the model structure on the category 
$\cC I\cat$ of $\cC I$-categories under slightly different assumptions,
namely when the underlying monoidal model category $\cC$ is simplicial and
has a monoidal fibrant replacement functor. Part (2) of the following 
Proposition is essentially the same as Lemma 3.6 of \cite{dundas}.

\begin{proposition}\label{prop-cat-mc}
Let $I$ be a fixed set and $\cC$ be a cofibrantly generated monoidal model 
category which satisfies the monoid axiom.  Assume as well that every object 
in $\cC$ is small. 
\begin{enumerate}
\item $\cC I\cat$ is a cofibrantly generated model category with
weak equivalences and fibrations defined pointwise.
\item  Every cofibration of $\cC I$-categories whose source is pointwise
cofibrant is a pointwise cofibration.  In particular, if the unit object 
$\mI_{\cC}$ is cofibrant in $\cC$, then
every cofibrant $\cC I$ category is pointwise cofibrant.
\end{enumerate}
\end{proposition}

\begin{proof} 
The generating (trivial) cofibrations are the image
of the free monoid functor $T_{\cC I}\co  \cC I\graph \to \cC I\cat$ 
applied to the sets of generating (trivial) cofibrations in $\cC I\graph$.
The proof of the first statement follows from~\cite[4.1(3), 6.2]{ss}.  The
third to last paragraph in the proof of~\cite[6.2]{ss} is the
only place which uses the symmetry of the monoidal product.  At that
point one is working in the underlying category, which is $\cC I\graph$ in  
our case.  Since both weak equivalences and pushouts in 
$\cC I\graph$ are determined pointwise,
one can just work pointwise in $\cC$.  Then the symmetry of the monoidal
product and the monoid axiom do hold by assumption.    

The proof of the second statement is essentially the same as in 
the last paragraph of \cite[6.2]{ss}, 
which treats the case of algebras, i.e., when $I$ is a singleton.
Note, this analysis does not require a symmetric monoidal product.
\end{proof}

As with Theorem~\ref{thm-WMQE gives equivalent algebras}
we use the following comparison of $\lambda$ with the structured
left adjoints $L^{\cC I}$ and $L_{\cO}$ for categories, respectively modules. 
In this paper we only apply Theorem \ref{thm-WMQE many generators} 
to {\em strong} monoidal Quillen pairs, namely in the next and final section.
In that case, the maps $\chi_{\cO}$ and $\chi_M$ are isomorphisms,
and so the proposition is redundant.   Elsewhere this statement
is used for weak monoidal Quillen pairs, though; see~\cite{sh-Q}.
Refer to Definition~\ref{def-create} for the notion of a `cell' object.
 
\begin{proposition}\label{prop-L-many}
Let  $(\!\xymatrix{\lambda\co \cD \ar@<-.4ex>[r]&\cC:R\ar@<-.4ex>[l]}\!)$ 
be a weak monoidal Quillen pair between monoidal model categories 
with cofibrant unit objects.
\begin{enumerate}
\item Suppose that the functor $R\co \cC I\cat\to\cD I\cat$  
has a left adjoint $L^{\cC I}$.
Then for every ${\cD}I$-cell category $\cO$,
the morphism of $\cC I$-graphs
\[ \chi_{\cO} \ \co \ \lambda \cO \ \to \ L^{\cC I}\cO \]  
which is adjoint to the underlying $\cD I$-graph morphism 
of the adjunction unit $\cO\to R(L^{\cC I}\cO)$
is a pointwise weak equivalence.
\item Suppose $\cO$ is a ${\cD}I$-cell category for which the
functor $R\co \Modr L^{\cC I}\cO \to \Modr \cO$ has a left adjoint $L_{\cO}$.
Then for every cell $\cO$-module $M$, the map 
\[ \chi_M \ \co \ \lambda M \ \to \ L_{\cO}M \]  
which is adjoint to the underlying $\cD$-morphism of the adjunction unit
$M\to R(L_{\cO}M)$ is an objectwise weak equivalence in $\cC$. 
\end{enumerate}
\end{proposition}

We omit the proof of Proposition \ref{prop-L-many}, since it is 
essentially identical to the proof of Proposition~\ref{prop-L}; 
the free monoid functor is replaced by the free $I$-category 
functor $T_{\cC I} \co\cC I\graph \to \cC I\cat$ and the
underlying objects are in $\cC I\graph$ and $\cD I\graph$ 
instead of $\cC$ and $\cD$. For modules there is another difference:
instead of one kind of free module, for every element 
$j\in I$ and every object $K$ of $\cC$, there is an $\cO$-module
$F_jK=K\tensor \cO(-,j)$ freely generated by $K$ at $j$.

Finally, we extend the results of Section~\ref{sec-main} to these enriched
categories.  

\begin{theorem}\label{thm-WMQE many generators}
Let $R\co \cC\to\cD$ be the right adjoint of a 
weak monoidal Quillen equivalence. Suppose that the unit objects in
$\cC$ and $\cD$ are cofibrant.

\begin{enumerate}
\item\label{ONE-many generators} 
Consider a cofibrant ${\cD}I$ category $\cO$ such that 
the forgetful functors create model structures
for modules over $\cO$ and modules over $L^{\cC I}\cO$.
Then the adjoint functor pair 
\[\xymatrix{L_{\cO}\ : \ \Modr \cO\ \ar@<-.4ex>[r]&
\ \Modr (L^{\cC I}\cO) \ : R\ar@<-.4ex>[l] }\]
is a Quillen equivalence.
\item\label{TWO-many generators}
Suppose that Quillen invariance holds for $I$-categories in 
$\cC$ and $\cD$. 
Then for any pointwise fibrant $\cC I$-category 
$\cA$ such that the forgetful functors 
create model structures for modules over $\cA$ and modules over $R\cA$,
the adjoint functor pair 
\[\xymatrix{L^{\cA}\ : \ \Modr R\cA\ \ar@<-.4ex>[r]&
\ \Modr \cA \ : R\ar@<-.4ex>[l] }\]
is a Quillen equivalence. 
If the right adjoint $R$ preserves all weak equivalences in $\cC$ 
and the forgetful functors create model structures for modules
over any monoid, then this holds for any $\cC I$-category $\cA$.
\end{enumerate}
\end{theorem}

The proof of Theorem~\ref{thm-WMQE many generators} is now
almost literally the same as the proof of parts 
\eqref{ONE} and \eqref{TWO}  
of Theorem \ref{thm-WMQE gives equivalent algebras}.
Whenever Proposition \ref{prop-L} is used in the proof of the latter
we now appeal to Proposition \ref{prop-L-many} instead.

\section{Symmetric monoidal categories of spectra}\label{sec-spectra} 

In this section we show that the Quillen equivalences between
the categories of rings, modules and algebras established in~\cite{mmss}
and~\cite{sch} can be extended to Quillen equivalences between
modules over `ring spectra with many objects' or `spectral categories'.  
This then shows that the classification results in~\cite{stable} 
can be translated to any one of the symmetric monoidal categories of spectra.

Comparison theorems between rings, modules and algebras based on
symmetric spectra 
over simplicial sets ($\spec_{sS}$) and topological spaces ($\spec_{Top}$), 
orthogonal spectra ($\cI\cS$), $\cW$-spaces ($\cW\cT$) 
and simplicial functors ($\cS\cF$)
can be found in~\cite{mmss}, Theorems 0.4 through 0.9 and 19.11.
Rings, modules and algebras based on $S$-modules $\cM_S$ are compared 
to their counterparts based on symmetric and orthogonal spectra 
in \cite[Thm.\ 5.1]{sch} respectively~\cite[Ch.~I]{mm}, 
Theorems 1.1 through 1.7.
See~\cite{sh-un} for an approach that unifies all of these comparison theorems.
These Quillen equivalences could also be deduced
via our general result in Theorem \ref{thm-WMQE gives equivalent algebras}. 
Moreover, in several of the categories which participate in the
diagram \eqref{spectra display} below, there are also Quillen equivalences for 
categories of {\em commutative} algebras
(with the exception of simplicial functors $\cS\cF$ 
and $\cW$-spaces $\cW\cT$); such results are out 
of the scope of the general methods in this present paper.
However, modules over `ring spectra with many objects' were not considered
in the above papers, and the point of this final section is to fill that gap.

The categories of spectra which we consider are all displayed 
in a commutative diagram of monoidal model categories and 
strong monoidal Quillen equivalences
\begin{equation}\label{spectra display}\xymatrix@C5mm@R8mm{ 
&\spec_{sS}\ar@<.3ex>[rr]^-{\mT} \ar@<-.3ex>[ddd]_{\mP} &&
\spec_{Top} \ar@<.3ex>[ll]^-{\mS}
\ar@<-.3ex>[dd]_{\mP} \ar@<-.3ex>[rr]_-{\Id} &&
(\spec_{Top})_+\ar@<.3ex>^-{\Lambda}[ddrr] \ar@<-.3ex>[ll]_-{\Id}
\ar@<-.3ex>[dd]_{\mP} \\
\\
&&& \cI\cS \ar@<-.3ex>[uu]_{\mU}\ar@<-.3ex>[rr]_-{\Id}\ar@<-.3ex>[d]_{\mP} &&
(\cI\cS)_+ \ar@<-.3ex>[uu]_{\mU} \ar@<-.3ex>[ll]_-{\Id} \ar@<.3ex>[rr]^{\mN}&&
\cM_S \ar@<.3ex>^-{\Phi}[uull] \ar@<.3ex>[ll]^{\mN^{\sharp}} \\
& \cS\cF\ar@<-.3ex>[uuu]_{\mU}\ar@<.3ex>[rr]^-{\mP\mT} &&
\cW\cT \ar@<-.3ex>[u]_{\mU}\ar@<.3ex>[ll]^-{\mS\mU}\\
\ar@<4ex>@{}[rr]|{\text{\ based on simplicial sets\ }} 
\ar@{(-)}[rrrr]|{\text{\ cofibrant unit object \ }} 
\ar@<-4ex>@{(-)}[rrrrrr]|{\text{\ diagram spectra \ }} & & \qquad
\ar@<4ex>@{(-)}[rrrrrr]|{\text{\ based on topological spaces\ }}& & \qquad
\ar@{(-)}[rrrr]|{\text{\ unit {\em not} cofibrant\ }} & &  \qquad & &  \qquad }
\medskip
\end{equation}
(where the left adjoints are on top and on the left).
There are five categories of {\em diagram spectra}:
symmetric spectra over simplicial sets $\spec_{sS}$~\cite{hss}, 
symmetric spectra over topological spaces $\spec_{Top}$~\cite{mmss}, 
simplicial functors $\cS\cF$~\cite{lydakis}, 
orthogonal spectra $\cI\cS$~\cite{mmss}, 
and $\cW$-spaces $\cW\cT$~\cite{mmss}.
The categories of topological symmetric spectra and orthogonal spectra
appear twice, with different model structures: the {\em stable}
model structure~\cite[3.4]{hss},~\cite[Sec.~9]{mmss} (without decoration) 
and the {\em positive stable} model structure~\cite[Sec.~14]{mmss}
(decorated with a subscript `+'). However, the stable and positive model
structures share the same class of weak equivalences.
The remaining category $\cM_S$ of $S$-modules~\cite{ekmm} is of a somewhat 
different flavor.
By $\mU$ we denote various `forgetful' or `underlying object' functors,
with the left adjoint `prolongation' functors $\mP$, which are all described
in the paper~\cite{mmss}.
Moreover, $\mS$ is the singular complex functor and $\mT$ is
geometric realization. 
The functors $\Lambda$ and $\Phi$ which relate symmetric spectra
to $S$-modules and their lifts $\mN$ and $\mN^{\sharp}$  to orthogonal spectra 
are defined and studied in~\cite{sch} respectively \cite[Ch.~I]{mm}.
See~\cite[4.7]{sh-un} for a unified approach to defining all of these
functors.

The following theorem is an application of Theorem \ref{a.1.2}
to these categories of spectra.

\begin{theorem} Let $\cC$ be any of the model categories
$\spec_{sS}$, $\spec_{Top}$, $ (\spec_{Top})_+$, $\cI\cS$, $ (\cI\cS)_+$, $
\cS\cF$, $\cW\cT$ or $\cM_S$. 
\begin{enumerate} 
\item The modules over any $\cC$-category 
inherit a model category structure in which the fibrations 
and weak equivalences are defined pointwise 
in the underlying category $\cC$.
\item If $\Psi\co \cO \to \cR$ 
is a pointwise weak equivalence of $\cC I$-categories, 
then restriction and extension of scalars along $\Psi$ form a $\cC$-Quillen
equivalence of the module categories.
\end{enumerate}
\end{theorem}
\begin{proof} (1)\qua All specified choices of monoidal model category $\cC$
are cofibrantly generated, 
see~\cite[12.1]{mmss}, ~\cite[3.4]{hss},~\cite[9.2]{lydakis} 
and \cite[VII.4]{ekmm}. 
For $\cM_S$ the argument for (1) follows just as for modules 
over a ring spectrum; see~\cite[VII.4.7]{ekmm}.  One could also verify
the pushout product, unit and monoid axioms directly. 
The unit axiom is automatic for the categories
$\spec_{sS}$, $\spec_{Top}$, $\cI\cS$, $\cS\cF$ and $\cW\cT$
where the unit object is cofibrant. In the positive model structures
$(\spec_{Top})_+$ and $(\cI\cS)_+$ the unit axiom holds 
since every positively cofibrant object is also stably cofibrant, 
and the respective unit objects are stably (but not positively) cofibrant.
Moreover, the pushout product and monoid axioms for the diagram spectra
other than $\cS \cF$ hold by~\cite[12.5, 12.6]{mmss},~\cite[5.3.8, 5.4.1]{hss}.
For the category of simplicial functors, these two axioms do not appear 
explicitly in Lydakis' paper~\cite{lydakis}, but we can argue as follows. 

By~\cite[12.3]{lydakis} the pushout product
of two cofibrations is a cofibration. To see that the pushout product
$i \boxprod j$ of a cofibration $i$ with a trivial cofibration $j$
is again a trivial cofibration, we argue indirectly and use the
Quillen equivalence $\mP\mT\co \cS\cF \to \cW\cT$ of~\cite[19.11]{mmss}.
By~\cite[3.5 (1)]{ss} or~\cite[4.2.5]{hovey-book} 
it suffices to check the pushout product of a generating cofibration 
with a generating trivial cofibration.
For the stable model structure of simplicial functors, 
Lydakis~\cite[9.1]{lydakis} uses generating cofibrations
and trivial cofibrations which all have cofibrant sources and targets.
Since the left Quillen functor $\mP\mT$ is also strong monoidal,
we have $\mP\mT(i \boxprod j) \iso (\mP\mT i) \boxprod (\mP\mT j)$
as morphisms of $\cW$-spaces. Since $\mP\mT i$ is a cofibration, 
$\mP\mT j$ is a trivial cofibration and the pushout product axiom holds in
$\cW\cT$, the pushout product $\mP\mT(i \boxprod j)$ is in particular
a stable equivalence of  $\cW$-spaces. 
As a left Quillen equivalence, $\mP\mT$ detects weak equivalences between
cofibrant objects, so $i\boxprod j$ is a stable equivalence of 
simplicial functors.

For the monoid axiom we consider a generating stable trivial cofibration $j$
from the set ${\bf SF}^g_{sac}$ defined in~\cite[9.1]{lydakis},
and we let $X$ be an arbitrary simplicial functor.
By~\cite[12.3]{lydakis}, $X\sm j$ is an injective morphism 
of simplicial functors; we claim that $X\sm j$ is also a stable equivalence.
To see this, we choose a cofibrant replacement $X^c\to X$;
then $X^c\sm j$ is a weak equivalence by the pushout product axiom.
By~\cite[12.6]{lydakis}, smashing with a cofibrant simplicial functor
preserves stable equivalences. Since the source and target of $j$ 
are cofibrant, $X\sm j$ is thus also a stable equivalence.
The class of injective stable equivalences of simplicial functors
is closed under cobase change and transfinite composition.
So every morphism in $({\bf SF}^g_{sac}\sm \cS\cF)\text{-cof}_{\text{reg}}$
is a stable equivalence, which implies the monoid axiom by~\cite[3.5 (2)]{ss}. 

In the categories $\spec_{sS}$ and $\cS\cF$ which are based on simplicial sets,
every object is small with respect to the whole category;
so the proof concludes by an application of Theorem \ref{a.1.2} (1).
In the other cases, which are based on topological spaces, the
{\em cofibration hypothesis} \cite[5.3]{mmss}, \cite[VII.4]{ekmm}
makes sure that the small object argument still applies and
the conclusion of Theorem \ref{a.1.2} is still valid.

Part (2) is proved by verifying the criterion 
of Theorem \ref{a.1.2} (2): for every cofibrant right $\cO$-module $N$, 
the induced map $N \tensor_{\cO} \cO \to N \tensor_{\cO} \cR$ 
is an objectwise weak equivalence.  
The special case of 
modules over a monoid, i.e., when the set $I$ has one element,
is treated in~\cite[12.7]{mmss}, \cite[5.4.4]{hss} and~\cite[III 3.8]{ekmm}.
Again for simplicial functors, this argument is not quite contained in 
\cite{lydakis}, but one can also verify the criterion 
for Quillen invariance as in~\cite[12.7]{mmss} 
using the fact that smashing with a cofibrant
object preserves stable weak equivalences by~\cite[12.6]{lydakis}.
The general case of modules over a category with more
than one object uses the same kind of cell induction as for
modules over a monoid; we omit the details.
\end{proof}

All the Quillen adjoint pairs appearing 
in the master diagram \eqref{spectra display}
are strong monoidal Quillen equivalences in the sense of
Definition \ref{def-WMQ equivalence}; so we can apply
Theorem \ref{thm-WMQE many generators}
to get Quillen equivalences for modules over `ring spectra with many objects';
this leads to the following proof of Corollary~\ref{cor-spectra}.
Special care has to be taken for the positive model structures
on the categories of symmetric and orthogonal spectra, and for
the category of $S$-modules since there the units of the smash product
are {\em not} cofibrant.

\begin{proof}[Proof of Corollary~\ref{cor-spectra}]
One of the main results of~\cite{stable}, Theorem 3.3.3, shows that
any cofibrantly generated, proper, simplicial stable model category
with a set of generators is Quillen equivalent to modules over a 
$\spec_{sS}$-category $\cO$ with one object for each generator. 
By Proposition~\ref{prop-cat-mc} 
we can choose a cofibrant replacement $\spec_{sS}$-cell category $\cO^c$ 
with a pointwise stable equivalence $q\co \cO^c\stackrel{\sim}{\to}\cO$.
Then $\Modr\cO$ and $\Modr\cO^c$ are Quillen equivalent.  

For comparisons which do not involve the category $\cM_S$ 
of $S$-modules, but only the left part of diagram \eqref{spectra display},
we consider the {\em stable model structures}.
In the five categories of diagram spectra, the unit of the smash product 
is cofibrant with respect to this stable model structure.
So various applications of part \eqref{ONE-many generators} 
of Theorem~\ref{thm-WMQE many generators} show 
that $\Modr \cO^c$ is Quillen equivalent to modules 
over the $\spec_{Top}$-category $\mT(\cO^c)$
(where $\mT$ is the geometric realization functor, applied levelwise),
to modules over the $\cI\cS$-category $\mP \mT(\cO^c)$,
to modules over the $\cW\cT$-category $\mP\mP\mT(\cO^c)$
(this composite $\mP\mP$ is just denoted by $\mP$ in~\cite{mmss}),
and to modules over the $\cS\cF$-category $\mP(\cO^c)$.

To connect to the world $\cM_S$ of $S$-modules we have to argue slightly
differently, since the unit $S$ in ${\mathcal M}_S$ is not cofibrant.
First we change model structures on the category of $\mT(\cO^c)$-modules
by viewing the identity functors as a Quillen equivalence between
the stable and positive model structures 
(which share the same class of weak equivalences); see the right-hand
part of the diagram below.

The last Quillen pair we consider compares modules 
over $\mT(\cO^c)$ and modules over the $\cM_S$-category
$\Lambda\mT(\cO^c)$. The right adjoint is given by $\Phi$,
together with restriction of scalars along the adjunction unit 
$\eta\co \mT(\cO^c)\to \Phi\Lambda\mT(\cO^c)$; the left adjoint is induced
by pointwise application of $\Lambda$.
The right adjoint $\Phi$ preserves all weak equivalences;
so to see that we have a Quillen equivalence we may show that
for every cofibrant $\mT(\cO^c)$-module $M$
the adjunction unit $M\to  \Phi\Lambda(M)$ is a pointwise stable equivalence.
Since $\cO^c$ is a $\spec_{sS}$-cell category, $\mT(\cO^c)$ is a  
$\spec_{Top}$-cell category, both times with respect to the
stable model structure on symmetric spectra.
The unit is cofibrant in the stable model structure of
symmetric spectra; hence $\mT(\cO^c)$ is pointwise cofibrant
in the stable model structure of symmetric spectra as well.

The {\em positive} cofibrations of symmetric spectra
are precisely those stable cofibration which are isomorphisms in
level 0. So every positively cofibrant $\mT(\cO^c)$-module $M$ is
also stably cofibrant. Since $\mT(\cO^c)$ itself is is pointwise 
stably cofibrant, so is $M$; thus the adjunction unit $M\to\Phi\Lambda(M)$ 
is a stable equivalence by \cite[Thm.\ 3.1]{sch}.

To sum up, we display all these Quillen equivalences in the
following diagram, where we also indicate the underlying monoidal
model categories:
\[\xymatrix@C9mm@R10mm{ 
\Modr \cO /\spec_{sS} \ar@<.3ex>[d]^-{q_*} & \\
\Modr \cO^c/ \spec_{sS}  \ar@<.3ex>[u]^-{q^*} 
\ar@<.3ex>[r]^-{\mT} \ar@<-.3ex>[dd]_{\mP} &
\Modr \mT(\cO^c)/ \spec_{Top} \ar@<.3ex>[l]^-{\mS}
\ar@<-.3ex>[d]_{\mP} \ar@<-.3ex>[r]_-{\Id} & 
\Modr \mT(\cO^c)/ (\spec_{Top})_+ 
\ar@<-.3ex>[d]_-{\Lambda} \ar@<-.3ex>[l]_-{\Id} \\
& \Modr \mP\mT(\cO^c)/\cI\cS \ar@<-.3ex>[u]_{\mU}\ar@<-.3ex>[d]_{\mP} &
\Modr \Lambda\mT(\cO^c)/ \cM_S\ar@<-.3ex>[u]_-{\Phi}\\
 \Modr \mP(\cO^c)/ \cS\cF\ar@<-.3ex>[uu]_{\mU}\ar@<.3ex>[r]^-{\mP\mT} &
\Modr \mP\mT\mP(\cO^c)/\cW\cT \ar@<-.3ex>[u]_{\mU}\ar@<.3ex>[l]^-{\mS\mU} }\]
\vspace{-20pt}
\end{proof}

\Addresses\recd

\begin{thebibliography}[EKMM]

\bibitem[BM]{Berger-Moerdijk}
C.~Berger and I.~Moerdijk, {\em Axiomatic homotopy theory for operads},
Preprint (2002). \verb!http://arXiv.org/abs/math.AT/0206094!

\bibitem[Bor94]{bor}
F.~Borceux, {\em Handbook of categorical algebra II, Categories and
structures}, Cambridge University Press, 1994.

\bibitem[BG76]{BG}
A.~K.~Bousfield and V.~K.~A.~M.~Gugenheim, {\em On ${\rm PL}$ de Rham
theory and rational homotopy type}, Mem. Amer. Math. Soc. {\bf 8}
(1976), no. 179, ix+94 pp.

\bibitem[Ca54]{cartan}
H.~Cartan, {\em Alg\`ebres d'Eilenberg-MacLane at homotopie},
S\'eminaire Henri Cartan, 1954-55

\bibitem[Do58]{dold}
A.~Dold, {\em Homology of symmetric products and other functors of complexes},
Ann. Math. {\bf 68} (1958), 54-80.

\bibitem[Du01]{dundas}
B.~I.~Dundas, {\em Localization of $V$-categories},
Theory Appl. Categ. {\bf 8} (2001), 284--312. 

\bibitem[Dw80]{Dwyer-operations}
W.~G.~Dwyer, {\em Homotopy operations for simplicial commutative algebras.}
Trans.\ Amer.\ Math.\ Soc.\ {\bf 260} (1980), 421--435. 

\bibitem[DK80]{DK}
W.~G. Dwyer and D.~M. Kan, {\em Function complexes in homotopical algebra},
Topology {\bf 19} (1980), 427--440.

\bibitem[DS95]{ds}
W.~G.~Dwyer and J.~Spalinski, {\em Homotopy theories and model categories}, 
Handbook of algebraic topology (Amsterdam), North-Holland, 
Amsterdam, 1995, pp.~73--126.

\bibitem[EM53]{EM-H(pi)}
S.~Eilenberg and S.~Mac Lane, {\em On the groups $H(\Pi,n)$, I},
Ann.\ of Math.\ (2) {\bf 58}, (1953), 55--106.

\bibitem[EKMM]{ekmm}
A.~D.~Elmendorf, I.~Kriz, M.~A.~Mandell, and J.~P.~May, 
{\em Rings, modules, and algebras in stable homotopy theory. 
{W}ith an appendix by M.~Cole},
Mathematical Surveys and Monographs, {\bf 47},  Amer.\ Math.\ Soc.,
Providence, RI, 1997, xii+249 pp.

\bibitem[GS]{GS}
J.~P.~C.~Greenlees and B.~Shipley, {\em Rational torus-equivariant
cohomology theories III: the Quillen equivalence}, in preparation.

\bibitem[Hov99]{hovey-book}
M.~Hovey, {\em Model categories}, Mathematical Surveys and Monographs,
{\bf 63}, Amer.\ Math.\ Soc., Providence, RI, 1999, xii+209 pp.

\bibitem[HSS]{hss}
M.~Hovey, B.~Shipley, and J.~Smith, {\em Symmetric spectra},
J.\ Amer.\ Math.\ Soc.\ {\bf 13} (2000), 149--208.

\bibitem[Jar97]{jardine}
J.~F.~Jardine,  {\em  A closed model structure 
for differential graded algebras},
Cyclic Cohomology and Noncommutative Geometry, 
Fields Institute Communications, {\bf 17}, AMS (1997), 55-58.

\bibitem[Lyd98]{lydakis}
M.~Lydakis, {\em Simplicial functors and stable homotopy theory},
Preprint (1998). \newline
\verb!http://hopf.math.purdue.edu/!

\bibitem[ML63]{ML-homology}
S.~Mac Lane, {\em Homology}, Grundlehren der math.\ Wissensch.\ {\bf 114},
Academic Press, Inc., Springer-Verlag,  1963 x+422 pp.

\bibitem[ML71]{ML-working}
S.~Mac Lane, {\em Categories for the working mathematician},
Graduate Texts in Math. {\bf 5}, 
Springer, New York-Berlin, 1971, ix+262 pp.

\bibitem[Man]{mandell}
M.~A.~Mandell, {\em Topological Andre-Quillen Cohomology and 
E-infinity Andre-Quillen Cohomology}, Adv. in Math., to appear.\newline
\verb!http://www.math.uchicago.edu/~mandell/!
 
\bibitem[MMSS]{mmss}
M.~A.~Mandell, J.~P.~May, S.~Schwede and B.~Shipley,
{\em Model categories of diagram spectra}, Proc.\ London Math.\ Soc.,
{\bf 82} (2001), 441-512.

\bibitem[MM02]{mm}
M.~A.~Mandell and J.~P.~May, {\em Equivariant orthogonal spectra and 
$S$-modules}, Memoirs Amer. Math. Soc., {\bf 159} (2002), no. 755, x+108 pp.

\bibitem[May67]{may}
J.~P.~May, {\em Simplicial objects in algebraic topology}, Chicago Lectures
in Mathematics, Chicago, 1967, viii+161pp.

\bibitem[Qui67]{Q}
D.~G. Quillen, {\em Homotopical algebra}, Lecture Notes in Mathematics,
{\bf 43}, Springer-Verlag, 1967.

\bibitem[Qui69]{Q2}
D.~G. Quillen, {\em Rational homotopy theory}, Ann. of Math.
{\bf 90} (1969), 204--265.

\bibitem[Ri03]{richter-DoldKan}
B.~Richter, {\em Symmetries of the Dold-Kan correspondence},
Math.\ Proc.\ Cambridge Phil.\  Soc.\ {\bf 134} (2003), 95--102.

\bibitem[Sch01]{sch}
S.~Schwede, {\em S-modules and symmetric spectra}, Math. Ann. {\bf 319}
(2001), 517-532. 

\bibitem[SS00]{ss}
S.~Schwede and B.~Shipley, {\em Algebras and modules 
in monoidal model categories}, Proc.\ London Math.\ Soc.\ {\bf 80} (2000),
491-511.

\bibitem[SS03]{stable}
S.~Schwede and B.~Shipley, 
{\em Stable model categories are categories of modules},
Topology, {\bf 42} (2003), 103-153.

\bibitem[Sh01]{sh-un}
B.~Shipley, {\em Monoidal uniqueness of stable homotopy theory},
Adv. in Math. {\bf 160} (2001), 217-240.

\bibitem[Sh02]{S-QT}
B.~Shipley, {\em An algebraic model for rational $S^1$-equivariant
stable homotopy theory}, Quart. J. of Math. {\bf 53} (2002),
87-110.

\bibitem[S]{sh-Q}
B.~Shipley, {\em $H\mZ$-algebra spectra are differential graded algebras},
Preprint (2002). \newline
\verb!http://www.math.purdue.edu/~bshipley/!

\bibitem[St]{stanley}
D.~Stanley, {\em Determining closed model category structures},
Preprint (1998). \newline
\verb!http://hopf.math.purdue.edu/!
\end{thebibliography}
\end{document}